\documentclass[11pt]{article}

\usepackage{amsmath}
\usepackage{amsfonts}
\usepackage{hyperref}
\usepackage{enumitem}
\usepackage[numbers]{natbib}
\usepackage{graphicx}
\usepackage{multirow}
\usepackage[table]{xcolor}
\usepackage{tikz}
\usetikzlibrary{shapes.geometric,arrows,positioning}
\tikzstyle{startstop} = [rectangle, rounded corners, minimum width=2cm, minimum height=0.6cm,text centered, draw=black, fill=red!20]
\tikzstyle{process}   = [rectangle, minimum width=2cm, minimum height=0.6cm, text centered, draw=black, fill=orange!20, align=center]
\tikzstyle{decision}  = [rectangle, rounded corners, aspect=2, minimum width=2cm, minimum height=0.6cm, text centered, draw=black, fill=blue!20, inner sep=5pt, align=center]
\tikzstyle{arrow}     = [thick,->,>=stealth]
\usepackage[ruled,vlined,linesnumbered]{algorithm2e}
\SetKwRepeat{Do}{do}{while}

\definecolor{lightgray}{gray}{0.9}

\DeclareMathOperator*{\argmin}{\mathrm{argmin}}

\begin{document}

\title{Models, constructive heuristics, and benchmark instances for
  the flexible job shop scheduling problem with sequencing flexibility
  and position-based learning effect\thanks{This work has been
    partially supported by the Brazilian agencies FAPESP (grants
    2013/07375-0, 2018/24293-0, and 2022/05803-3) and CNPq (grants
    311536/2020-4 and 302073/2022-1).}}

\author{
    K. A. G. Araujo\thanks{Department of Applied Mathematics, Institute of
    Mathematics and Statistics, University of S\~ao Paulo, Rua do
    Mat\~ao, 1010, Cidade Universit\'aria, 05508-090, S\~ao Paulo, SP,
    Brazil. e-mail: kennedy94@ime.usp.br}
 \and 
    E. G. Birgin\thanks{Department of Computer Science, Institute of
    Mathematics and Statistics, University of S\~ao Paulo, Rua do
    Mat\~ao, 1010, Cidade Universit\'aria, 05508-090, S\~ao Paulo, SP,
    Brazil. e-mail: egbirgin@ime.usp.br}
 \and
    D. P. Ronconi\thanks{Department of Production Engineering, Polytechnic 
    School, University of S\~ao Paulo, Av. Prof. Luciano Gualberto, 1380, 
    Cidade Universit\'aria, 05508-010 São Paulo, SP, Brazil. e-mail: 
    dronconi@usp.br}
}

\date{February 26, 2024}

\maketitle

\begin{abstract}
This paper addresses the flexible job shop scheduling problem with
sequencing flexibility and position-based learning effect. In this
variant of the flexible job shop scheduling problem, precedence
constraints of the operations constituting a job are given by an
arbitrary directed acyclic graph, in opposition to the classical case
in which a total order is imposed. Additionally, it is assumed that
the processing time of an operation in a machine is subject to a
learning process such that the larger the position of the operation in
the machine, the faster the operation is processed. Mixed integer
programming and constraint programming models are presented and
compared in the present work. In addition, constructive heuristics are
introduced to provide an initial solution to the models' solvers. Sets
of benchmark instances are also introduced. The problem considered
corresponds to modern problems of great relevance in the printing
industry. The models and instances presented are intended to support
the development of new heuristic and metaheuristics methods for this
problem.\\

\end{abstract}

\section{Introduction}

In this work, we consider the flexible job shop (FJS) scheduling
problem with sequence flexibility and position-based learning
effect. The problem is NP-hard, as it has the job shop scheduling
problem, known to be NP-hard~\cite{Garey1976}, as a particular
case. The sequencing flexibility feature refers to the fact that the
precedence constraints imposed on the operations of a job are given by
an arbitrary directed acyclic graph instead of the usual linear order
enforced in the FJS scheduling. Many real-life scheduling problems fit
into this scope, such as, for example, problems in the printing
industry~\cite{birgin2014milp,Lunardi2020,Lunardi2021}, mold
manufacturing industry~\cite{Gan2002}, and glass
industry~\cite{AlvarezValdes2005}. In~\cite{Lunardi2021}, the FJS with
sequencing flexibility and several additional features, such as, for
example, resumable operations, periods of unavailability of the
machines, sequence-dependent setup times, partial overlapping between
operations with precedence constraints, and fixed operations, was
addressed. However, learning effects have not yet been taken into
account in the literature regarding FJS with sequencing
flexibility. The present work is devoted to this problem and presents
the first step towards its effective and efficient resolution. Mixed
integer linear programming (MILP) and constraint programming (CP)
models and a relatively large set of instances are introduced for the
purpose of producing a benchmark test set. Since, leaving aside very
small instances, commercial exact solvers alone are hardly even able
to find a feasible solution, constructive heuristics are proposed with
the goal of enhancing their performance. Overall, this work takes a
first step towards solving the proposed problem and provides a solid
and robust basis for the future development of more sophisticated
heuristic and metaheuristic methods.

In classical scheduling problems, the processing time of a given
operation on a given machine is a fixed input parameter. However,
there are real-life situations in which the manufacturing process
involves repetitive manual tasks and the worker undergoes a learning
process that results in a reduction of the execution time of his/her
task. For instance, workers can get proficient at performing
assemblies quickly, or more confident and skillful in manipulating
hardware, software, and/or raw materials. If we consider that the
reduction in the operation processing time is related to the number of
times the worker has already performed the operation, we are dealing
with a position-based learning effect. The pioneering works in
applying the concept of learning effect to scheduling problems
are~\cite{Biskup1999,Cheng2000,Gupta1988}. Surveys on the subject can
be found in~\cite{Azzouz2017,Biskup2008,Janiak2011}.

A brief literature review in chronological order of the FJS with
sequencing flexibility follows. In~\cite{Gao2006}, the problem was
addressed by considering non-fixed intervals of machine unavailability
for preventive maintenance. A multi-objective MILP formulation and a
hybrid multi-objective genetic algorithm were
proposed. In~\cite{Gan2002} and~\cite{Kim2003}, process flexibility
was also taken into account, which means that the same result can be
obtained with different operation sequences. In~\cite{Gan2002} a
branch-and-bound algorithm was proposed, while~\cite{Kim2003}
considered a symbiotic evolutionary
algorithm. In~\cite{VitalSoto2020}, where a MILP formulation was
presented to minimize the weighted tardiness, a hybrid bacterial
foraging optimization algorithm was developed. Furthermore, the
algorithm was enhanced by a local search method based on the
manipulation of critical operations. A research that addresses the FJS
with sequencing flexibility and sequence-dependent setup times can be
found in~\cite{Cao2021}. The authors proposed a knowledge-based cuckoo
search algorithm associated with a reinforcement learning strategy for
self-adjustment. In~\cite{Kasapidis2021}, MILP and CP models and a
hybrid evolutionary algorithm with local search mechanisms were
introduced. A variation in which the processing of each operation
requires multiple resources was considered in~\cite{Kasapidis2023}, in
which models are presented and some properties of the problems are
analyzed.

In~\cite{AlvarezValdes2005}, an application in the glass industry was
described and a heuristic approach combining local search and priority
rules was proposed to minimize the total cost related to final product
completion times. Other industrial environments, such as the printing
industry, have also been modeled as an FJS with sequencing
flexibility. Regarding this particular application, \cite{Vilcot2008}
suggested a bi-objective genetic algorithm based on NSGA II to solve
the problem. In~\cite{birgin2014milp}, a MILP model and a constructive
heuristic were presented, while~\cite{Birgin2015} introduced a list
scheduling algorithm and its extension to a beam search method. In
\cite{Lunardi2020,Lunardi2021}, formulations using constrained
programming and mixed integer linear programming were established, as
well as trajectory and population metaheuristics were introduced.
In~\cite{Yu2017}, it was considered the flight deck scheduling problem,
which is an FJS with sequencing flexibility with additional
constraints that state that some operations must be completed before
others. The problem was described through its MILP formulation and
instances were solved with a differential evolution type method.
In~\cite{AndradePineda2020}, the scheduling of repair orders and
worker assignment in an automotive repair shop was analyzed. The main
scheduling problem is a two-resource FJS with sequencing
flexibility. The problem was modeled by extending the formulation
introduced in~\cite{birgin2014milp} and an iterated greedy heuristic
was presented.

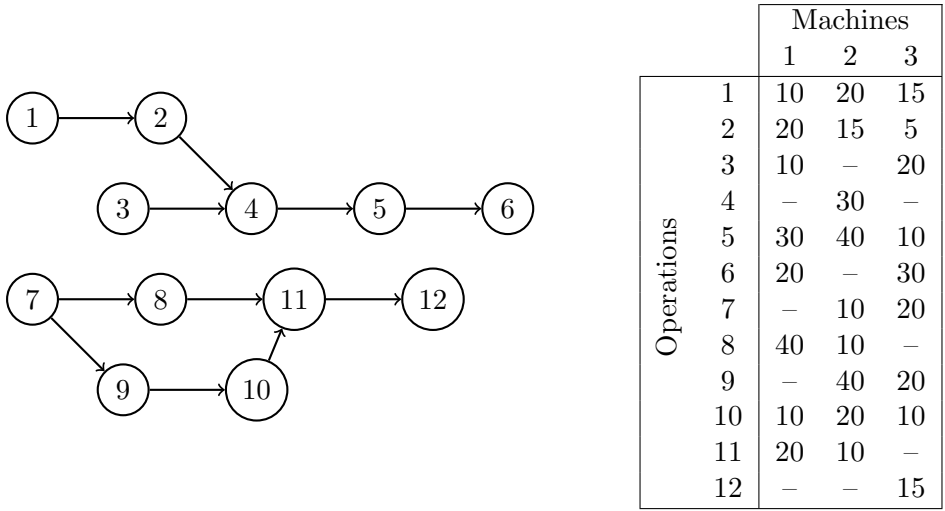
\begin{figure}[ht!]
\begin{center}
\begin{tabular}{cc}
\begin{minipage}{8cm}
\begin{tikzpicture}[node distance={10mm}, thick, main/.style = {draw, circle}] 
\node[main] (1) at (0,0) {1};
\node[main] (2) [right =of  1] {2};
\node[main] (3) [below right =of  1] {3};
\node[main] (4) [below right =of 2] {4};
\node[main] (5) [right =of 4] {5};
\node[main] (6) [right =of 5] {6};
\node[main] (7) [below left = of 3]  {7};
\node[main] (8) [right =of 7] {8};
\node[main] (9) [below right =of 7] {9};
\node[main] (10) [right =of  9] {10};
\node[main] (11) [right =of  8] {11};
\node[main] (12) [right =of  11] {12};
\draw[->] (1)--(2);
\draw[->] (2)--(4);
\draw[->] (3)--(4);
\draw[->] (4)--(5);
\draw[->] (5)--(6);
\draw[->] (7)--(8);
\draw[->] (7)--(9);
\draw[->] (8)--(11);
\draw[->] (9)--(10);
\draw[->] (10)--(11);
\draw[->] (11)--(12);
\end{tikzpicture}
\end{minipage}
&
\begin{minipage}{5cm}
\begin{tabular}{|cc|ccc|}
\cline{3-5}
\multicolumn{2}{c|}{} & \multicolumn{3}{c|}{Machines}\\
\multicolumn{2}{c|}{} & 1 & 2 & 3\\
\hline
\multirow{12}{*}{\rotatebox{90}{Operations}}
&  1 & 10 & 20 & 15\\
&  2 & 20 & 15 &  5\\
&  3 & 10 & -- & 20\\
&  4 & -- & 30 & --\\
&  5 & 30 & 40 & 10\\
&  6 & 20 & -- & 30\\
&  7 & -- & 10 & 20\\
&  8 & 40 & 10 & --\\
&  9 & -- & 40 & 20\\
& 10 & 10 & 20 & 10\\
& 11 & 20 & 10 & --\\
& 12 & -- & -- & 15\\
\hline
\end{tabular}
\end{minipage}
\end{tabular}
\end{center}
\caption{On the left, representation of operations' precedence
  constraints by a directed acyclic graph $D=(\mathcal{O},\widehat
  A)$, where $\mathcal{O} = \{ 1, 2, \dots, 12 \}$ represents the set
  of operations and $\widehat A$ is the set of arcs that represent the
  precedence constraints. On the right, standard processing times of
  the twelve operations on each of the three machines. In the table
  cells, ``--'' means that the machine cannot process the operation.}
\label{fig1}
\end{figure}

Let us consider an illustrative example of the FJS with sequencing
flexibility and position-based learning effect. The example has~12
operations and~3 machines. Figure~\ref{fig1} shows the precedence
relationships between the operations and the standard processing times
of each operation on each machine. The small table with the standard
processing times shows that in the FJS scheduling problem each
operation can be processed by one or more machines (situation known as
routing flexibility), in opposition to the job shop (JS) scheduling
problem in which each operation can be processed by only one
machine. The concept of job is implicitly defined by the directed
acyclic graph (DAG) and corresponds to a set of operations that have
some dependency relationship between them. The figure makes it clear
that, in the FJS scheduling problem with sequence flexibility, the
dependencies of the operations of a job are given by an arbitrary
acyclic directed graph as opposed to the total order of the FJS
scheduling problem. In this example there are two jobs, one with 6
operations (numbered from 1 to 6) and the other also with 6 operations
(numbered from 7 to 12).

\begin{figure}[ht!]
\centering
\begin{tikzpicture}[node distance={12mm}, thick, main/.style = {draw, circle}] 
\node[main] (12) {$s$}; 
\node[main] (14) [above right =of 12,label={[orange]90:15}] {1};
\node[main] (1) [right =of  14,label={[orange]90:5}] {2};
\node[main] (2) [right =of  12,label={[cyan]-90:10}] {3};
\node[main] (3) [below right =of 1,label={[violet]90:30}] {4};
\node[main] (4) [right =of 3,label={[orange]90:10}] {5};
\node[main] (5) [right =of 4,,label={[cyan]90:20}] {6};
\node[main] (6) [below right= of 12,label={[violet]-90:10}]  {7};
\node[main] (7) [right =of 6,label={[cyan]-90:40}] {8};
\node[main] (8) [below right =of 6,label={[orange]-90:20}] {9};
\node[main] (9) [right =of  8,,label={[orange]-90:10}] {10};
\node[main] (10) [below =of  4,label={[violet]-90:10}] {11};
\node[main] (11) [below =of  5,label={[orange]90:15}] {12};
\node[main] (13) [above right =of 11] {$t$};
\draw[->,yellow,ultra thick] (12) -- (14) -- (1) -- (3) -- (4) -- (5) -- (13);
\draw[->,purple] (12) -- (14);
\draw[->,purple] (12) -- (2);
\draw[->,purple] (12) -- (6);
\draw[->] (14)--(1);
\draw[->] (1)--(3);
\draw[->] (2)--(3);
\draw[->] (3)--(4);
\draw[->] (4)--(5);
\draw[->] (6)--(7);
\draw[->] (6)--(8);
\draw[->] (7)--(10);
\draw[->] (8)--(9);
\draw[->] (9)--(10);
\draw[->] (10)--(11);
\draw[->,purple] (11)--(13);
\draw[->,purple] (5)--(13);
\draw[->,cyan,dashed] (7)--(5);
\draw[->,violet,dashed] (6)--(3);
\draw[->,violet,dashed] (3)--(10);
\draw[->,orange,dashed] (1) to [bend right] (8);
\draw[->,orange,dashed] (8) to [bend left] (9);
\draw[->,orange,dashed] (9)--(4);
\draw[->,orange,dashed] (4)--(11);
\draw[->,orange,dashed] (14) to [bend left] (1);
\draw[->,cyan,dashed] (2)--(7);
\end{tikzpicture} 
\caption{Representation of an optimal solution to the instance in
  Figure~\ref{fig1}.}
\label{fig2}
\end{figure}
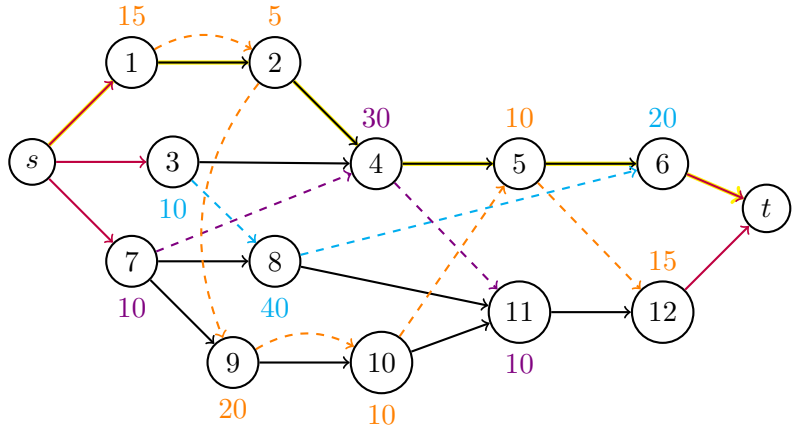

\begin{figure}[ht!]
\centering
\includegraphics[scale=0.6]{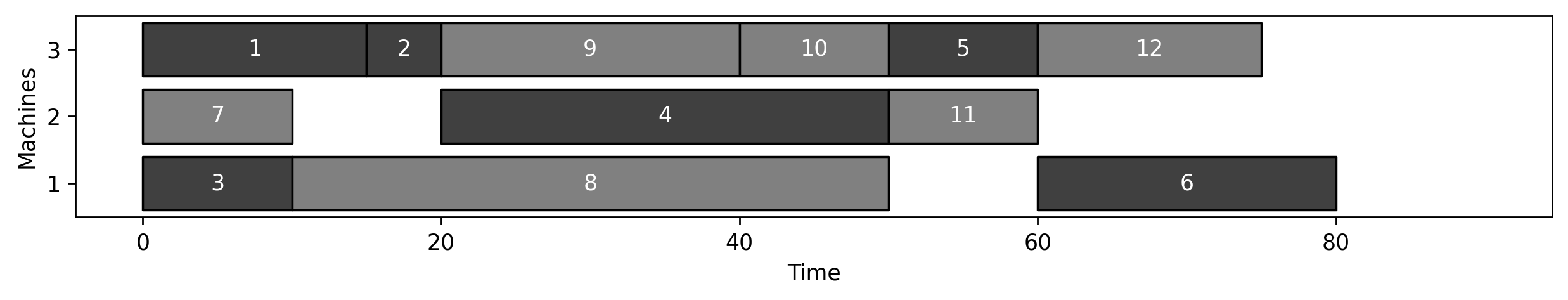}
\caption{Gantt chart representation of the optimal solution shown in
  Figure~\ref{fig2} to the instance of Figure~\ref{fig1}.}
\label{fig3}
\end{figure}

A feasible solution to the instance of Figure~\ref{fig1} can be
illustrated by a DAG in which a source node~$s$ and a target node~$t$
are added to the DAG that represents the precedence constraints; see
Figure~\ref{fig2}. Arcs from~$s$ to all operations with no
predecessors and from all nodes without successor to~$t$ must also be
added. (For the readers that can see the figure in colors, those arcs
are painted purple.) In addition, dashed arcs represent the sequence
(list) in which operations are processed by each machine. (For the
readers that can see the figure in colors, blue corresponds to
machine~1, violet corresponds to machine~2, and orange corresponds to
machine~3.) Colored figures at the top or bottom of the operations
correspond to their processing times. In this directed graph, that we
name $G=(V,A)$ from now on\footnote{Note that $V=\mathcal{O} \cup
\{s,t\}$ and $A$ is composed by the given arcs in $\widehat A$ that
represent the precedence relations among operations, plus the
mentioned arcs related to the new nodes $s$ and $t$, plus the
mentioned machine arcs.}, the longest path from $s$ to $t$ corresponds
to the makespan. In this example, the longest path, with value~80,
corresponds to the path $s, 1, 2, 4, 5, 6, t$. (Highlighted yellow in
the figure for those readers who can see the figure in color.) The
depicted feasible solution corresponds to an optimal solution. Figure
\ref{fig3} shows the Gantt chart representation of the solution. Note
that standard times were used, i.e.\ in this example the learning
effect was not considered at all.

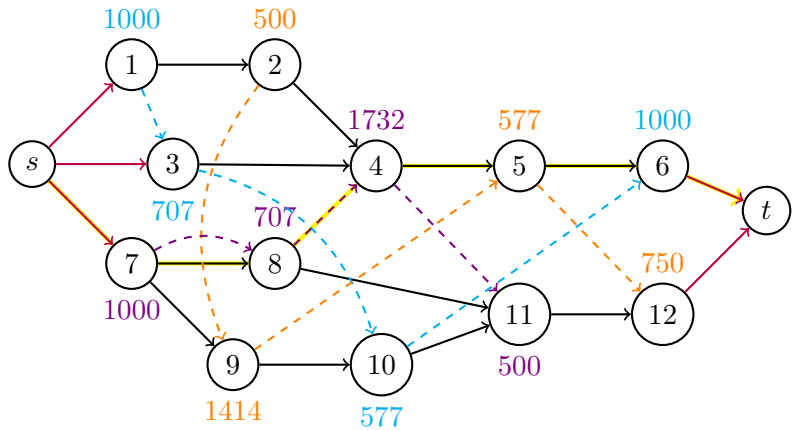
\begin{figure}[ht!]
\centering
\begin{tikzpicture}[node distance={12mm}, thick, main/.style = {draw, circle}] 
\node[main] (13) {$s$}; 
\node[main] (1) [above right =of 13,label={[cyan]90:1000}] {1};
\node[main] (2) [right =of  1,label={[orange]90:500}] {2};
\node[main] (3) [right =of  13,label={[cyan]-90:707}] {3};
\node[main] (4) [below right =of 2,label={[violet]90:1732}] {4};
\node[main] (5) [right =of 4,label={[orange]90:577}] {5};
\node[main] (6) [right =of 5,,label={[cyan]90:1000}] {6};
\node[main] (7) [below right= of 13,label={[violet]-90:1000}]  {7};
\node[main] (8) [right =of 7,label={[violet]90:707}] {8};
\node[main] (9) [below right =of 7,label={[orange]-90:1414}] {9};
\node[main] (10) [right =of  9,label={[cyan]-90:577}] {10};
\node[main] (11) [below =of  5,label={[violet]-90:500}] {11};
\node[main] (12) [below =of  6,label={[orange]90:750}] {12};
\node[main] (14) [above right =of 12] {$t$};
\draw[->,yellow,ultra thick] (13) -- (7) -- (8) -- (4) -- (5) -- (6) -- (14);
\draw[->,purple] (13) -- (1);
\draw[->,purple] (13) -- (3);
\draw[->,purple] (13) -- (7);
\draw[->] (1)--(2);
\draw[->] (2)--(4);
\draw[->] (3)--(4);
\draw[->] (4)--(5);
\draw[->] (5)--(6);
\draw[->] (7)--(8);
\draw[->] (7)--(9);
\draw[->] (8)--(11);
\draw[->] (9)--(10);
\draw[->] (10)--(11);
\draw[->] (11)--(12);
\draw[->,purple] (12)--(14);
\draw[->,purple] (6)--(14);
\draw[->,cyan,dashed] (1)--(3);
\draw[->,cyan,dashed] (3) to [bend left] (10);
\draw[->,cyan,dashed] (10)-- (6);
\draw[->,violet,dashed] (8)--(4);
\draw[->,violet,dashed] (4)--(11);
\draw[->,orange,dashed] (2) to [bend right] (9);
\draw[->,orange,dashed] (9)--(5);
\draw[->,orange,dashed] (5)--(12);
\draw[->,violet,dashed] (7) to [bend left] (8);
\end{tikzpicture} 
\caption{Representation of an optimal solution to the instance in
  Figure~\ref{fig1} in the presence of learning effect.}
\label{fig4}
\end{figure}

\begin{figure}[ht!]
\centering
\includegraphics[scale=0.6]{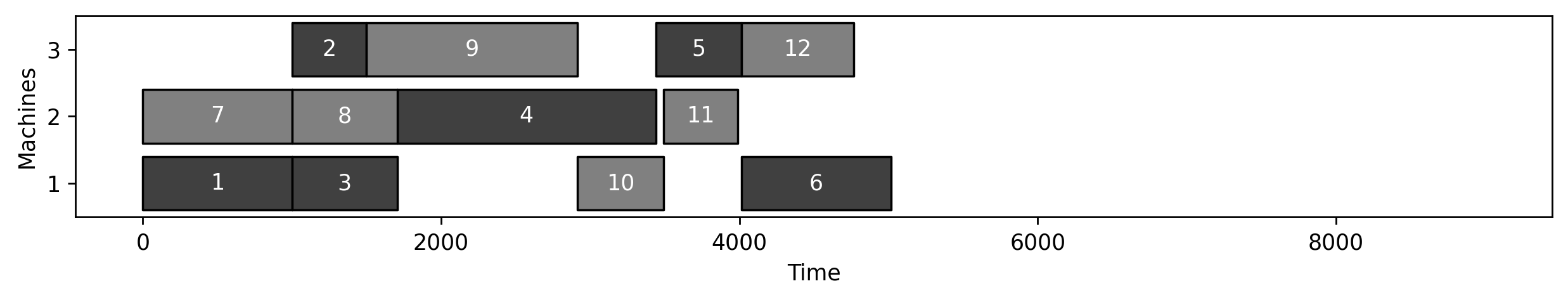}
\caption{Gantt chart representation of the optimal solution shown in
  Figure~\ref{fig4} to the instance of Figure~\ref{fig1} in the
  presence of learning effect.}
\label{fig5}
\end{figure}

In the present work, we consider that when an operation~$i$ is
assigned to a machine~$k$ and it is the $r$th operation to be
processed by the machine, the standard processing time $p_{ik}$ is
affected by a learning effect and becomes $\psi_{\alpha}(p_{ik},r)$,
where $\alpha>0$ is the learning rate. Following~\cite{Biskup1999}, we
might consider $\psi_{\alpha}(p,r) := p / r^{\alpha}$ in which, the
larger the value of $\alpha$, the faster the learning. However, in the
present work, a constraint programming model of the problem will be
introduced and instances solved with a commercial solver that requires
processing times to assume integer values only. Thus, we consider
$\psi_{\alpha}(p,r) := \lfloor 100 \, p / r^{\alpha} + 1/2
\rfloor$. Multiplying by one hundred, adding~$0.5$, and taking the
floor corresponds to changing the unit of measures (from seconds to
milliseconds, for example) and rounding to the closest integer
value. Figures~\ref{fig4} and~\ref{fig5} illustrate the graph
representation and the Gantt chart of an optimal solution in which the
learning function $\psi_{\alpha}$ with $\alpha=0.5$ is considered. It
is worth noting that a different schedule is found, whose makespan,
given by the critical path $s, 7, 8, 4, 5, 6, t$, is equal to 50.16 in
the original units of time (5016 in the new one).

The rest of this paper is organized as follows. In
Section~\ref{models}, we introduce the integer linear programming and
constraint programming models. In Section~\ref{heuristics}, we outline
the proposed constructive heuristics. In Section~\ref{instances}, we
report the introduced instances. Numerical experiments with the
constructive heuristics and exact commercial solvers are reported in
Section~\ref{experiments}. Conclusions and lines of future research
are presented in the concluding section.

\section{Mixed integer and constraint programming models} \label{models}

In this section, we present mixed-integer linear programming (MILP)
and constraint programming (CP) formulations for the FJS scheduling
problem with sequencing flexibility and position-based learning
effect.

\subsection{Mixed-integer linear programming model}

The adoption of position-based decision variables serves as the
fundamental approach for modeling problems involving position-based
learning effects, as it enables a more natural expression of
constraints related to the change in processing times. The proposed
MILP model is derived from~\cite{birgin2014milp} and is built upon the
model introduced in~\cite{wagner1959integer} that considers
position-based decision variables; see
also~\cite{wilson1989alternative}. Notation below simplifies the
presentation of the model.
\begin{description}[topsep=2pt,itemsep=0pt]
\item[Sets:] $\phantom{}$
\begin{description}[topsep=0pt,itemsep=0pt]
\item[$\mathcal{O}$:] set of operations,
\item[$\mathcal{F}$:] set of machines,
\item[$\mathcal{O}_k$:] set of operations that can be processed by machine $k$,
\item[$\mathcal{F}_i$:] set of machines that can process operation $i$,
\item[$\widehat A$:] set of directed arcs that represent operations' precedence constraints.
\end{description}

\item[Parameters:] $\phantom{}$
\begin{description}[topsep=0pt,itemsep=0pt]
\item[$p_{ik}$:] standard processing time of operation $i$ in machine $k$,
\item[$\psi_{\alpha}(p_{ik},r)$:] position-based learning function.
\end{description}

\item[Decision variables:] $\phantom{}$
\begin{description}[topsep=0pt,itemsep=0pt]
\item[$x_{ikr}$:] 1 if operation $i$ is the $r$th operation to be processed in machine $k$; 0, otherwise,
\item[$s_i$:] starting time of operation $i$,
\item[$h_{kr}$:] starting time of the $r$th operation to be processed in machine $k$,
\item[$p'_{i}$:] actual processing time of operation $i$ (considering the learning effect).
\end{description}
\end{description}
The MILP model follows.
\begin{align}
\mbox{Minimize} \quad C_{\max} \quad & \label{m1_fo}
\\
\mbox{subject to} \quad & \nonumber \\
\sum_{k \in \mathcal{F}_i} \sum_{r = 1}^{|\mathcal{O}_k|} x_{ikr}  = 1, & \quad  i \in \mathcal{O}, \label{m1_1}\\
\sum_{i \in \mathcal{O}_k} x_{ikr} \leq 1, & \quad k \in \mathcal{F},r = 1,\dots,|\mathcal{O}_k|, \label{m1_2}\\
\sum_{i \in \mathcal{O}_k} x_{i,k,r+1} \leq \sum_{i \in \mathcal{O}_k} x_{ikr}, & \quad k \in \mathcal{F}, r = 1, \dots, |\mathcal{O}_k|-1, \label{m1_3}\\
p'_i = \sum_{k \in \mathcal{F}_i} \sum_{r = 1}^{|\mathcal{O}_k|}  \psi_{\alpha}(p_{ik},r) \, x_{ikr}, & \quad   i \in \mathcal{O}, \label{m2_1}\\
h_{kr} + \sum_{i \in \mathcal{O}_k}  \psi_{\alpha}(p_{ik},r) \, x_{ikr} \leq h_{k,r+1}, & \quad k \in \mathcal{F}, r = 1, \dots, |\mathcal{O}_k| - 1, \label{m2_2}\\
h_{kr} + \sum_{i \in \mathcal{O}_k}  \psi_{\alpha}(p_{ik},r) \, x_{ikr} \leq C_{\max}, & \quad k \in \mathcal{F}, r = |\mathcal{O}_k|, \label{m2_3}\\
s_i + p'_i \leq s_j, & \quad \forall (i,j) \in \widehat A, \label{m1_7}\\
s_i + p'_i - \left(2 - x_{ikr} - \sum_{t = r +1}^{|\mathcal{O}_k|} x_{jkt}\right) M \leq s_j,
& \quad \begin{array}{c}\forall i,j \in \{\mathcal{O} \;|\; i \neq j\}, \forall k \in \mathcal{F}_i \cap \mathcal{F}_j,\\ r = 1,\dots, |\mathcal{O}_k| - 1, \end{array} \label{m1_8}\\
h_{kr} - M \, (1 - x_{ikr} )\leq s_i, & \quad i \in \mathcal{O},  k \in \mathcal{F}_i,r = 1, \dots, |\mathcal{O}_k|, \label{m1_9}\\
s_i - M \, (1 - x_{ikr} )\leq h_{kr}, & \quad i\in \mathcal{O}, k \in \mathcal{F}_i,r = 1, \dots, |\mathcal{O}_k|, \label{m1_10}\\
s_i \geq 0, & \quad  i \in \mathcal{O}, \label{m1_11}\\
h_{kr} \geq 0, & \quad k \in \mathcal{F}_i, r = 1, \dots, |\mathcal{O}_k|, \label{m1_11b}\\
x_{ikr} \in \{0,1\}, & \quad  i \in \mathcal{O},k \in \mathcal{F}_i, r = 1,\dots,|\mathcal{O}_k|. \label{m1_12}
\end{align}

Objective function~\eqref{m1_fo} represents the minimization of the
makespan. Constraints~\eqref{m1_1} state that each operation must be
processed by exactly one machine and occupy exactly one position on
that machine.  Constraints~\eqref{m1_2} state that each position of
each machine can be associated with at most one
operation. Constraints~\eqref{m1_3} say that a position of a machine
can be occupied by an operation only if the previous position is also
occupied by another operation. Constraints~\eqref{m2_1} define the
actual processing time of each operation~$i$ in order to simplify the
presentation of the model. Constraints~\eqref{m2_2} avoid the
overlapping of operations assigned to the same
machine. Constraints~\eqref{m2_3} say that the makespan must be
greater than or equal to the completion time of all
operations. Combining these constraints with the minimization of the
objective function~\eqref{m1_fo} makes the makespan to coincide with
the completion time of the last operation to be completed.
Constraints~\eqref{m1_7} impose the given precedence constraints
between operations by saying that if an operation~$i$ precedes an
operation~$j$ then~$j$ cannot be started before~$i$ is completed.
Constraints~\eqref{m1_8} state that, if two operations~$i$ and~$j$
were assigned to the same machine~$k$ and operation~$i$ precedes
operation $j$, then~$j$ cannot be started before~$i$ is completed. If
operation~$i$ is assigned to position~$r$ of machine~$k$ (i.e.\ if
$x_{ikr}=1$) then constraints~\eqref{m1_9} and~\eqref{m1_10} force
$h_{kr}$ (the starting time of $r$th operation to be processed by
machine~$k$) to be equal to $s_i$ (the starting time of
operation~$i$). Constraints~\eqref{m1_11}, \eqref{m1_11b}, and~\eqref{m1_12} determine
the decision variables' domains. $M$ is a ``sufficiently large''
number and may assume the value $\sum_{i\in \mathcal{O}} \sum_{k\in
  {\cal F}} p_{ik}$. Function $\psi_{\alpha}(p_{ik},r)$ is the given
function that represents the learning effect and computes the actual
processing time of operation~$i$ if it is assigned to position~$r$ of
machine~$k$. This function has the learning rate~$\alpha \geq 0$ as a
parameter and, in the present work, as already mentioned in the
introduction, is given by
\[
\psi_{\alpha}(p,r) = \left\lfloor 100 \; r^{-\alpha} \; p + \frac{1}{2} \right\rfloor.
\]

\subsection{Constraint Programming model}

Constraint Programming (CP) is a potent methodology widely employed
for solving scheduling problems in academic and industrial
literature. CP Optimizer~\cite{Laborie2018}, being an optimization commercial solver
rooted in CP, incorporates specialized concepts of constraints and
variables, significantly facilitating the modeling process for
scheduling problems. In this section, we introduce a CP model
specifically designed for its utilization in connection with CP
Optimizer. The syntax of CP Optimizer is defined as it arises within
the formulation. The model follows below.
\begin{align}
\mbox{Minimize}\quad \max_{i \in \mathcal{O}} \text{endOf}(o_i) \label{CP_fo} 
\\
\mbox{subject to} \quad & \nonumber \\
\text{endBeforeStart}(o_i,o_j), & \quad (i,j) \in \widehat A, \label{CP_1}\\
\text{alternative}\left(o_i, [a_{ikr}]_{k \in \mathcal{F}_i, r = 1, \dots,|\mathcal{O}_k|}\right), & \quad i \in \mathcal{O},\label{CP_2}\\
\text{noOverlap}\left([a_{ikr}]_{i \in \mathcal{O}_k, r = 1, \dots,|\mathcal{O}_k| }\right), & \quad k \in \mathcal{F},\label{CP_3}\\
\text{endBeforeStart}(a_{ikr},a_{jk,r+1}), &
\quad \begin{array}{c} i,j \in \mathcal{O}, \, k \in \mathcal{F}_i \cap \mathcal{F}_j, \\ r = 1,\dots, |\mathcal{O}_k| - 1, \end{array} \label{forca2}\\
\text{or}(\text{[presenceOf}(a_{ik,r+1})]_{i \in \mathcal{O}_k}) \implies \text{or}(\text{[presenceOf}(a_{ikr}])_{i \in \mathcal{O}_k}), & \quad k \in \mathcal{F}, r = 1, \dots, |\mathcal{O}_k|-1, \label{CP_4}\\
\text{interval } o_i, & \quad i \in \mathcal{O},\label{CP_5}\\
\text{interval } a_{ikr}, \text{ opt}, \text{ size } = \psi_{\alpha}(p_{ik},r), & \quad i \in \mathcal{O}, k \in \mathcal{F}_i, r = 1, \dots, |\mathcal{O}_k|. \label{CP_6} 
\end{align}

Interval decision variables of the problem are described in
\eqref{CP_5} and \eqref{CP_6}. In \eqref{CP_5}, an interval variable
$o_i$ for each operation $i$ is defined. In \eqref{CP_6}, an
\textit{optional} interval variable $a_{ikr}$ is defined for each
possible assignment of operation $i$ to a machine $k \in {\cal F}_i$
at positions $r = 1, \dots, |\mathcal{O}_k|$. \textit{Optional} means
that the interval variable may exist or not; and the remaining of the
constraint says that, in case it exists, its size must be given by
$\psi_{\alpha}(p_{ik},r)$. Constraints \eqref{CP_2} state that each
operation $i$ must be allocated to exactly one machine $k \in
\mathcal{F}_i$ in exactly one position $r$, that is, one and only one
interval variable $a_{ikr}$ must be present and the selected interval
$a_{ikr}$ must start and end at the same instants as interval
$o_i$. The objective function \eqref{CP_fo} is to minimize the
makespan, given by the maximum end value of all the operations
represented by the interval variables $o_i$. Precedence constraints
between operations are posted as endBeforeStart constraints between
interval variables in constraints \eqref{CP_1}. Constraints
\eqref{CP_3} state that, for each machine $k$, at most one operation
can be assigned to each position and operations assigned to different
positions cannot overlap. Constraints~\eqref{forca2} say that, for any
machine and any position, the operation assigned to that position must
finish before the operation in the next position starts. Constraints
\eqref{CP_4} force that the empty positions of machine $k$ are the
last ones, i.e., that operations are processed in the first positions
of the machines without empty positions among them.


\section{Constructive Heuristics} \label{heuristics}

In this section, we propose two constructive heuristics for the FJS
scheduling problem with sequencing flexibility and position-based
learning effect. Constructive heuristics are algorithms that build a
feasible solution from scratch by iteratively selecting and sequencing
one operation at a time. The two proposed constructive heuristics are
based on the earliest starting time (EST) rule~\cite{birgin2014milp}
and the earliest completion time (ECT) rule~\cite{Leung2005}. The goal
is to use them to provide an initial feasible solution to the exact
solver that will be used to attempt to solve a set of test instances.

Algorithm~\ref{alg1} presents the constructive heuristic based on EST
rule. In the algorithm, $r^{\mathrm{op}}_v$ refers to the ready time
of operation~$v$, $w_v$ refers to its actual processing time and $c_v$
to its completion time. On the side of the machines,
$r^{\mathrm{mac}}_k$ represents the instant in which machine~$k$ is
released and~$g_k$ represents its first free position, which is the
one that would be occupied if an operation were assigned to it (both
quantities refer to the partial scheduling being constructed). $f_v$
will indicate to which machine operation $v$ was assigned and each
machine $k$ will have an ordered list $Q_k$ with the sequence of
operations to be processed. After the initializations (lines~2 to~5)
comes the main loop, which is executed as long as there are still
unscheduled operations. Among the not scheduled ones, the ready time
is calculated for all those that already have all the preceding
operations scheduled (lines~7 to~9). In line~10, observing the ready
times of the operations and machines, the smallest instant $r_{\min}$
in which an operation could be scheduled is calculated and the set $E$
of operation/machine pairs that could start at that instant $r_{\min}$
is constructed. As it was observed in~\cite{birgin2014milp}, $|E|$ can
be quite large and experience shows that a tie-breaking rule can
significantly improve the method's performance. Thus, in line~11,
among all the operation/machine pairs in $E$, taking into account the
learning effect, the pair $(\hat v,\hat k)$ with the shortest
processing time is chosen. In line 12, $w_{\hat v}$, $f_{\hat v}$, and
$c_{\hat v}$ are defined and the ready time $r^{\mathrm{mac}}_{\hat
  k}$ and the free position~$g_{\hat k}$ of machine~$\hat k$ are
updated. In lines~13 and~14 the corresponding machine arc is inserted
in the graph~$G$ (the arc must not be inserted if operation~$\hat v$
is the first one of machine~$\hat k$). Finally, the list of operations
assigned to machine~$\hat k$ is updated and the scheduled operation is
removed from the set of operations not yet scheduled. After all the
operations have been scheduled, the critical path in $G$ is calculated
(line~16) to determine the makespan value~$C_{\max}$. This is done
with Algorithm~\ref{alg2}. Initializations in lines 2 to 5 of
Algorithm~\ref{alg1} have complexity $O(|\mathcal{O}| +
|\mathcal{\widehat A}| + |\mathcal{F}|)$. Within the main loop (lines
6 to 15), lines 7--9 have complexity $O(|\mathcal{A}|)=O(|\widehat
A|+|\mathcal{O}|+|\mathcal{F}|)$, lines 10--11 have complexity
$O(|\mathcal{O}| + \sum_{i \in \mathcal{O}} |F_i|)$, and line 12 has
complexity $O(|\mathcal{O}|)$. Since the main loop is executed
$|\mathcal{O}|$ times, its total complexity is $O(|\mathcal{O}|
(|\widehat A|+|\mathcal{F}|+\sum_{i \in \mathcal{O}}
|\mathcal{F}_i|))$. As the complexity of the main loop is larger that
the complexity of the initialization as well as the complexity of line
16 (see below), then the complexity of Algorithm~\ref{alg1} is given
by the complexity of its main loop. It is important to note that
$\gamma = \sum_{i \in \mathcal{O}} |\mathcal{F}_i|$ is between
$|\mathcal{O}|$ and $|\mathcal{O}| |\mathcal{F}|$, but we prefer to
keep the complexity expressed as a function of $\gamma$ because
$\gamma$ is a measure of the size of the input that depends on the
sequencing flexibility of the instance under consideration. It is also
important to note that the complexity of the algorithm depends on the
routing flexibility of the operations, i.e., it depends on the number
of dependency relations in $\widehat A$. Therefore, it is important to
represent an instance in such a way that $\widehat A$ corresponds to a
transitive reduction of the precedences' digraph.


\begin{algorithm}[!ht]
\caption{Computes a solution graph $G=(V,A)$, $f$, $Q$, and $w$ by
  using EST dispatching rule. Then, in $G$, it computes the largest
  path $\mathcal{P}$ from $s$ to $t$ and its length~$C_{\max}$.}
\label{alg1}
\KwIn{${\cal O}$, $\mathcal F$, $p$, $\widehat A$}
\KwOut{$f$, $w$, $Q$, $G = (V, A)$, $\mathcal{U}$, $\mathcal{P}$, $C_{\max}$, $\tau$}
\SetKwBlock{Begin}{function}{end function}
\DontPrintSemicolon
\Begin({EST(}${\cal O}$, ${\cal F}$, $p$, $\widehat A$, $f$, $w$, $Q$, $G$, $\mathcal{U}$, $\mathcal{P}$, $C_{\max}$, $\tau${)}){
    
Set $A \gets \widehat A \cup \{ (s,j) \;|\; (\cdot,j) \not\in \widehat A \} \cup \{
(i,t) \;|\; (i,\cdot) \not\in \widehat A \}$ and define $V := \mathcal{O}
\cup \{s,t\}$ and $G = (V,A)$.\;

Set $r^{\mathrm{op}}_v \gets +\infty$ for all $v \in V$
and define $r^{\mathrm{op}}_s := 0$, $w_s := w_t := 0$, and
$c_s:=0$. \;

Set $r^{\mathrm{mac}}_k \gets 0$ and $g_k \gets 1$ for all $k \in
\mathcal{F}$.\;

Initialize $\Pi \gets V \setminus\{s,t\}$ as the set of
non-scheduled operations, and $Q_k$ as an empty list for all $k
\in \mathcal{F}$.\;

\While{$\Pi \neq \emptyset$}{

    \For{$v \in \Pi$}{
        \If{$\Pi \cap \{i \;|\:(i,v) \in A \}  = \emptyset$}{
            $r^{\mathrm{op}}_v \gets \max \{ c_i \mid i \in V\setminus \Pi$ such that $(i,v) \in A \}$\;}} 
        
    Set $r_{\min} = \min \{ \max(
    r^{\mathrm{op}}_v,r^{\mathrm{mac}}_k) \;|\; v \in \Pi, k \in
    \mathcal F_v \}$ and let $E$ be the set of pairs $(v,k)$ with
    $v \in \Pi$ and $k \in \mathcal F_v$ such that $\max(
    r^{\mathrm{op}}_v,r^{\mathrm{mac}}_k) = r_{\min}$.\;
 
    $(\hat v,\hat k) \gets \argmin \{ r_{\min} +
    \psi_{\alpha}(p_{v, k}, g_k) \;|\; (v,k) \in E\}$.\;
        
    Define $w_{\hat v} := \psi_{\alpha}(p_{\hat v,\hat k}, g_{\hat
      k})$, $f_{\hat v} := \hat k$ and $c_{\hat v} :=
    \max(r^{\mathrm{op}}_{\hat v}, r^{\mathrm{mac}}_{\hat k}) +
    w_{\hat v}$, and set $r^{\mathrm{mac}}_{\hat k} \gets c_{\hat
      v}$ and $g_{\hat k} \gets g_{\hat k} + 1 $.\;
 
    \If{$|Q_{\hat k}| \neq 0$}{
        Let $Q_{\hat k} = i_1, \dots, i_{|Q_{\hat k}|}$. Set $A \gets A \cup \{ (i_{|Q_{\hat k}|},\hat v) \} $.\;}
    
    Insert $\hat v$ at the end of $Q_{\hat k}$ and set $\Pi \gets \Pi \setminus \{\hat v\}$. \;
    }

CriticalPath($\mathcal{F}$, $f$, $w$, $Q$, $G$, $\mathcal{U}$,
$\mathcal{P}$, $C_{\max}$, $\tau$).
}
\end{algorithm}

The critical path in the directed graph~$G=(V,A)$ can be computed with
an adaptation~\cite[\S 22.2]{Cormen:2009:IA:1614191} of the
Bellman-Ford algorithm; see Algorithm \ref{alg2}. In addition to the
critical path~$\mathcal{P}$, the algorithm returns a topological
order~$\mathcal{U}$ of the vertices of~$G$ and a vector~$\tau$ of
dimension~$|\mathcal{F}|$. The vector~$\tau$ stores, in element
$\tau_k$, the largest position in the list $Q_k$ (list of operations
assigned to machine~$k$) that contains an operation in the critical
path. These two elements are not used in the context of the present
work, but they are useful information for developing neighborhoods in
local search algorithms. Algorithm~\ref{alg2} has worst case time
complexity $O(|V| + |A|)$. Since $|V|=O(|\mathcal{O}|)$ and $|A| =
O(|\widehat A| + |\mathcal{O}| + |\mathcal{F}|)$, this complexity
translates into $O(|\mathcal{O}| + |\widehat A| + |\mathcal{F}|)$.
The topological order of Algorithm~\ref{alg2} is calculated with the
help of Algorithm~\ref{alg2b}, which implements a depth-first
search. Algorithm~\ref{alg2b} computes, optionally, for each $v \in
V$, the set $R_v^{\leftarrow}$ formed by the vertices~$w$ such that
there exists in $G$ a path from~$w$ to~$v$, i.e.\ the vertices that
can reach~$v$. Algorithm~\ref{alg2b} is recursive and its complexity
is $O(|V|+|\widehat A|)$. The sets $R_v^{\leftarrow}$ are not used in
the context of the present work, but, again, they contain valuable
information for the development of local search strategies.


\begin{algorithm}[ht!]
\caption{Computes a critical path $\mathcal{P}$ and its length $\xi$
  for a given graph $G = (V, A)$. In addition, if $\tau$ is present as
  an input parameter, determines the last critical operation in each
  machine.}
\label{alg2}
\KwIn{$\mathcal{F}$, $f$, $w$, $Q$, $G = (V,A)$, $\tau$}
\KwOut{$\mathcal{U}$, $\mathcal{P}$, $\xi$, $\tau$}
\SetKwBlock{Begin}{function}{end function}
\DontPrintSemicolon
\Begin($\text{CriticalPath}{(}\mathcal{F}$, $f$, $w$, $Q$, $G$, $\mathcal{U}$, $\mathcal{P}$, $\xi$, $\tau{)}$){
    
Initialize $d_i \leftarrow -\infty$ for all $i \in
V\setminus\{s\}$ and define $d_s := 0$ and $\pi_s := 0$.\;


Initialize $\mathcal{V} \leftarrow \emptyset$ and $\mathcal{U}$ as an empty list and compute in $\mathcal{U}$ a topological sort of the vertices in $V$, by calling TopologicalSort+($G$, $\mathcal{U}$, $s$, $\mathcal{V}$).\;

\For{$\ell=1,\dots,|V|$}{
    Let $i$ be the $\ell$-th operation in the topological order given by $\mathcal{U}$.\;
    \For{$j\; \mathrm{such} \; \mathrm{that} \;(i,j) \in A$}{
        \If{$d_j < d_i + w_i$}{
            $d_j \gets d_i + w_i$ and $\pi_j \gets i$.\;}}}

$\xi := d_t$\;

Initialize $i \gets \pi_t$, $\mathcal{P} \leftarrow \emptyset$,
and $\tau_k \leftarrow 0$ for all $k \in \mathcal{F}$.\;

\Do{$i \neq s$}{
  \If{$\tau_{f_i} = 0$}{
      Let $Q_{f_i}$ be given by the sequence $i_1,
      \dots,i_{\ell-1},i,i_{\ell+1},\dots,i_{|Q_{f_i}|}$. Define
      $\tau_{f_i} := \ell$.}
  
  $\mathcal{P} \gets \mathcal{P} \cup \{i\}$ and $i \gets \pi_i$.}
}
\end{algorithm}

\begin{algorithm}[!ht]
\caption{Computes a topological sort $\mathcal{U}$ of the vertices of $G = (V, A)$. In addition, if~$\mathcal{R}_v^{\leftarrow}$ is present as an input parameter, computes the set $\mathcal{R}_v^{\leftarrow}$ of vertices that reaches $v$ in $G = (V, A)$.} 
\label{alg2b}
\KwIn{$G = (V,A)$, $\mathcal{U}$, $v$, $\mathcal{V}$, $\mathcal{R}_v^{\leftarrow}$}
\KwOut{$\mathcal{U}$, $\mathcal{V}$, $\mathcal{R}_v^{\leftarrow}$}
\DontPrintSemicolon
\SetKwBlock{Begin}{function}{end function}
\Begin(TopologicalSort+{(}$G$, $\mathcal{U}$, $v$, $\mathcal{V}$, $\mathcal{R}_v^{\leftarrow}${)})
{
Set $\mathcal{V} \gets \mathcal{V} \cup \{v\}$.\;
\For{$j \; \mathrm{such} \; \mathrm{that} \; (v,j) \in A$}{
    \If{$j \notin \mathcal{V}$}{
        TopologicalSort+($G$, $\mathcal{U}$,  $j$, $\mathcal{V}$, $\mathcal{R}_v^{\leftarrow}$)\;
        }
    \If{$i \not\in \mathcal{R}_v^{\leftarrow} \;\mathrm{and} \;  j \in \mathcal{R}_v^{\leftarrow}$}{
        set $\mathcal{R}_v^{\leftarrow} \gets \mathcal{R}_v^{\leftarrow} \cup \{ v\}$.\;}
    }
    Insert $i$ at the beginning of $\mathcal{U}$.\;
}
\end{algorithm}

Algorithm~\ref{alg3} presents the constructive heuristic based on the
ECT rule. The algorithm is very similar to Algorithm~\ref{alg1},
except for one detail. In the constructive heuristic based on EST, we
first compute the instant~$r_{\min}$ which is the earliest instant at
which an unscheduled operation could be initiated. All
operation/machine pairs that could start at that instant are
considered and the pair with the shortest processing time is
selected. But since they would all start at instant~$r_{\min}$, saying
that the pair with the shortest processing time is chosen is the same
as saying that the pair that ends earliest is selected. This is the
idea that is taken to the extreme in the constructive heuristic based
on the ECT rule: without limiting the choice to the operation/machine
pairs that could start as early as possible, the operation/machine
pair that will finish earliest is chosen, even if the processing of
the operation does not start as early as possible. The worst case time
complexity of Algorithm~\ref{alg3} is the same as that of
Algorithm~\ref{alg1}.

\begin{algorithm}[!ht]
\caption{Computes the solution graph $G=(V,A)$, $f$, $Q$, and $w$ by
  using ECT dispatching rule. Then, in $G$, it computes the largest
  path $\mathcal{P}$ from $s$ to $t$ and its length~$C_{\max}$.}
\label{alg3}
\KwIn{${\cal O}$, $\mathcal F$, $p$, $\widehat A$}
\KwOut{$f$, $w$, $Q$, $G = (V, A)$, $\mathcal{U}$, $\mathcal{P}$, $C_{\max}$, $\tau$}
\SetKwBlock{Begin}{function}{end function}
\DontPrintSemicolon
\Begin(ECT{(}${\cal O}$, ${\cal F}$, $p$, $\widehat A$, $f$, $w$, $Q$, $G$, $\mathcal{U}$, $\mathcal{P}$, $C_{\max}$, $\tau${)}){
    
Set $A \gets \widehat A \cup \{ (s,j) \;|\; (\cdot,j) \not\in \widehat A \} \cup \{
(i,t) \;|\; (i,\cdot) \not\in \widehat A \}$ and define $V := \mathcal{O}
\cup \{s,t\}$ and $G = (V,A)$.\;

Set $r^{\mathrm{op}}_v \gets +\infty$ and define
$r^{\mathrm{op}}_s := 0$, $w_s := w_t := 0$, and $c_s:=0$. \;

Set $r^{\mathrm{mac}}_k \gets 0$ and $g_k \gets 1$ for all $k \in
\mathcal{F}$.\;

Initialize $\Pi \gets V \setminus\{s,t\}$ as the set of
non-scheduled operations, and $Q_k$ as an empty list for all $k
\in \mathcal{F}$.\;

\While{$\Pi \neq \emptyset$}{

    \For{$v \in \Pi$}{
        \If{$\Pi \cap \{i \;|\:(i,v) \in A \}  = \emptyset$}{
            $r^{\mathrm{op}}_v \gets \max \{ c_i \mid i \in V\setminus \Pi$ such that $(i,v) \in A \}$\;}} 

    $(\hat v,\hat k) \gets \argmin \{ \max(
    r^{\mathrm{op}}_v,r^{\mathrm{mac}}_k) + \psi_{\alpha}(p_{v, k}, g_k)
    \;|\; v \in \Pi, k \in \mathcal F_v\}$.\;
        
    Define $w_{\hat v} := \psi_{\alpha}(p_{\hat v,\hat k}, g_{\hat
      k})$, $f_{\hat v} := \hat k$ and $c_{\hat v} :=
    \max(r^{\mathrm{op}}_{\hat v}, r^{\mathrm{mac}}_{\hat k}) +
    w_{\hat v}$ and set $r^{\mathrm{mac}}_{\hat k} \gets c_{\hat
      v}$ and $g_{\hat k} \gets g_{\hat k} + 1 $.\;
    
    \If{$|Q_{\hat k}| \neq 0$}{
      Let $Q_{\hat k} = i_1, \dots, i_{|Q_{\hat k}|}$. Set $A
      \gets A \cup \{ (i_{|Q_{\hat k}|},\hat v) \} $.\;}
    
    Insert $\hat v$ at the end of $Q_{\hat k}$ and set $\Pi \gets
    \Pi \setminus \{\hat v\}$.\;
}

CriticalPath($\mathcal{F}$, $f$, $w$, $Q$, $G$, $\mathcal{U}$,
$\mathcal{P}$, $C_{\max}$, $\tau$).
}
\end{algorithm}

The EST-based heuristic gives priority to those pairs
operations/machines that can start the earliest. At the beginning of
the construction, this corresponds, roughly, to giving priority to all
the first operations of each job, which are operations that have no
precedents (operations~1, 3 and~7 in the example of
Figure~\ref{fig1}). Still, due to the intention of scheduling operations as
early as possible, it is possible that preference is given to empty
machines, building solutions that use several machines. By rapidly
scheduling the first operations of each job, more operations come to
have their precedents scheduled, increasing the number of
possibilities (search space) in future iterations of the method. On
the other hand, the heuristic based on the ECT rule chooses the
operation/machine pairs that terminate the earliest, regardless of
whether they are the ones that can start the earliest or not. Such a
strategy can limit the number of operation/machine pairs available in
future iterations, reducing the search space of the method. Moreover,
the choice for the operation/machine pair that can finish earliest,
combined with the learning effect, leads the method to schedule
operations on machines that already have several operations assigned
to them, since the higher the position in the machine, the shortest de
processing (reduced by the positioned-based learning effect). This
leads to the construction of solutions in which not all machines are
used. Depending on the learning rate~$\alpha$ considered and the
density of the DAG of precedences of the instance at hand, one
heuristic may be better than the other.

\section{Benchmark instances} \label{instances}

Tractability of the introduced models and performance of the proposed
constructive heuristics will be evaluated with the~50 large-sized
instances proposed in~\cite{birgin2014milp}, that were introduced for
the FJS scheduling problem with sequence flexibility but without
learning effect. In addition to these large-sized instances, a new set
with~60 smaller instances, using the generator described
in~\cite{birgin2014milp}, was generated. Instances whose name starts
with ``Y'' correspond to instances in which DAGs that represent the
operations' precedences are Y-shaped (like the DAG in the top of
Figure~\ref{fig1}); while instances whose name starts with ``DA''
correspond to instances in which DAGs that represent the operations
precedences are arbitrary DAGs (like the DAG in the bottom of
Figure~\ref{fig1}). The former will be called instances of Y-type and
the latter will be called instances of DA-type from now on.

For a given instance, we define measures $\omega_1$ and $\omega_2$ of
the sequencing flexibility and the routing flexibility,
respectively. Both measures are between~0 and~1 and, the larger their
values, the larger the flexibility their represent. Moreover, the
larger the flexibility, the larger the search space and, in
consequence, the harder the instance. Let $D^+ = (\mathcal{O},\widehat
A^+)$ be the transitive closure of the precedences DAG $D =
(\mathcal{O},\widehat A)$. Let $n_{\kappa}$ be the number of
operations of job~$\kappa$ and let $a_{\kappa}$ the number of
arcs in the $D^+$ among them. Then, $a^{\min}_\kappa \leq a_{\kappa}
\leq a^{\max}_{\kappa}$, where $a^{\min}(n_{\kappa})=n_{\kappa}-1$ and
$a^{\max}(n_{\kappa})=n_{\kappa}(n_{\kappa}-1)/2$. Therefore
\[
\omega_1^{\kappa} := 1 - \frac{a_{\kappa} -
  a^{\min}(n_{\kappa})}{a^{\max}(n_{\kappa} ) - a^{\min}(n_{\kappa})},
\]
is such that $\omega_1^{\kappa} \in [0,1]$ represents the sequencing
flexibility of the operations of job~$\kappa$. The arithmetic mean among
all~$n$ jobs of the instance, given by
\[
\omega_1 = \frac{1}{n} \sum_{\kappa=1}^n \omega_1^{\kappa},
\]
is also between 0 and 1 and represent the degree of the instance
sequencing flexibility. The larger $\omega_1$, the larger the
search space and, therefore, the more difficult the instance. (Of
course, any other type of average could be used in the definition of
$\omega_1$.)

In a similar way, we define the routing flexibility measure
$\omega_2$. It is easy to see that $|\mathcal{O}| \leq \sum_{i \in
  \mathcal{O}} |\mathcal{F}_i| = \sum_{k \in \mathcal{F}} |\mathcal{O}_k|
\leq |\mathcal{O}| |\mathcal{F}|$. Therefore,
\[
\omega_2 = \frac{ \sum_{i \in \mathcal{O}} |\mathcal{F}_i| -
  |\mathcal{O}|}{|\mathcal{O}| |\mathcal{F}| - |\mathcal{O}| },
\]
is between 0 and 1. The closer to~1, the larger the search space and,
therefore, the harder the instance.

Tables~\ref{tab1} and~\ref{tab2} provide a comprehensive overview of
the instances main characteristics. In the tables, $|\mathcal{O}|$ is
the number of operations, $|\mathcal{F}|$ is the number of machines,
$n$ is the number of jobs (connected components in the precedence
constraints DAG $D=(\mathcal{O},\widehat A)$), $|\widehat A|$ is the
number of precedence constraints, $\sum_{i \in \mathcal{O}}
|\mathcal{F}_i|$ is the number of ``true'' entries in the logical
matrix of dimension $|\mathcal{O}| \times |\mathcal{F}|$ whose $(i,k)$
says whether operation~$i$ can be processed by machine~$k$, and
$\omega_1$ and $\omega_2$ are the measures of the sequencing and
routing flexibility described in the previous paragraph. In addition,
for the MILP and the CP models, the tables show the number of
variables and constraints.

In the small-sized instances of Table~\ref{tab1}, the number of binary
variables of the MILP models and the number of interval variables of
the CP models go up to almost $1{,}000$; while in both models the
number of constraints goes up to $13{,}000$. On the other hand, in the
large-sized instances of Table~\ref{tab2}, the number of binary
variables of the MILP models and the number of interval variables of
the CP models go up to almost $73{,}000$; while in both models the
number of constraints goes up to $4{,}000{,}000$. Moreover, for each
instance, the number of binary variables in its MILP model is very
similar to the number of interval variables in its CP model and the
two models also have a very similar number of constraints.

\begin{table}[ht!]
\centering
\resizebox{0.8\textwidth}{!}{
\rowcolors{3}{}{lightgray}
\begin{tabular}{ccccccccrrrrr} 
\hline
\multicolumn{8}{c}{Main instance characteristics} &
\multicolumn{3}{c}{MILP formulation} &
\multicolumn{2}{c}{CP Optimizer formulation} \\ 
\hline
Instance & $|\mathcal{F}|$ & $|\mathcal{O}|$ & $n$ & $|\widehat A|$ & $\sum_{i \in \mathcal{O}} \mathcal{F}_i$ & $\omega_1$ & $\omega_2$ &
\begin{tabular}[c]{@{}r@{}}\#binary\\ variables\end{tabular} &
\begin{tabular}[c]{@{}r@{}}\#continous\\ variables\end{tabular} & \#constraints &
\begin{tabular}[c]{@{}r@{}}\#interval\\ variables\end{tabular} & \#constraints \\ 
\hline
miniDAFJS01 & 5 & 14 & 2 & 14 & 44 & 0.54 & 0.54 & 392 & 78 & 3,754 & 406 & 3,564 \\
miniDAFJS02 & 5 & 11 & 2 & 9 & 35 & 0.57 & 0.55 & 251 & 63 & 2,018 & 262 & 1,867 \\
miniDAFJS03 & 5 & 10 & 2 & 12 & 31 & 0.50 & 0.53 & 197 & 57 & 1,435 & 207 & 1,301 \\
miniDAFJS04 & 5 & 9 & 2 & 9 & 29 & 0.25 & 0.56 & 175 & 53 & 1,234 & 184 & 1,109 \\
miniDAFJS05 & 5 & 15 & 2 & 15 & 38 & 0.39 & 0.38 & 310 & 74 & 2,881 & 325 & 2,714 \\
miniDAFJS06 & 5 & 14 & 2 & 14 & 46 & 0.20 & 0.57 & 448 & 80 & 4,742 & 462 & 4,544 \\
miniDAFJS07 & 5 & 11 & 2 & 10 & 34 & 0.23 & 0.52 & 240 & 62 & 1,924 & 251 & 1,777 \\
miniDAFJS08 & 5 & 9 & 2 & 8 & 24 & 0.42 & 0.42 & 122 & 48 & 770 & 131 & 665 \\
miniDAFJS09 & 5 & 15 & 2 & 13 & 47 & 0.63 & 0.53 & 465 & 83 & 5,018 & 480 & 4,815 \\
miniDAFJS10 & 5 & 11 & 2 & 10 & 32 & 0.28 & 0.48 & 208 & 60 & 1,530 & 219 & 1,391 \\
miniDAFJS11 & 5 & 18 & 2 & 16 & 43 & 0.45 & 0.35 & 409 & 85 & 4,371 & 427 & 4,181 \\
miniDAFJS12 & 5 & 12 & 2 & 12 & 38 & 0.55 & 0.54 & 298 & 68 & 2,602 & 310 & 2,438 \\
miniDAFJS13 & 5 & 10 & 2 & 8 & 27 & 0.58 & 0.43 & 149 & 53 & 979 & 159 & 861 \\
miniDAFJS14 & 5 & 14 & 2 & 12 & 44 & 0.50 & 0.54 & 394 & 78 & 3,800 & 408 & 3,610 \\
miniDAFJS15 & 5 & 11 & 2 & 11 & 35 & 0.40 & 0.55 & 261 & 63 & 2,242 & 272 & 2,091 \\
miniDAFJS16 & 5 & 13 & 2 & 12 & 42 & 0.18 & 0.56 & 358 & 74 & 3,302 & 371 & 3,121 \\
miniDAFJS17 & 5 & 11 & 2 & 11 & 30 & 0.13 & 0.43 & 194 & 58 & 1,491 & 205 & 1,360 \\
miniDAFJS18 & 5 & 11 & 2 & 10 & 34 & 0.38 & 0.52 & 234 & 62 & 1,798 & 245 & 1,651 \\
miniDAFJS19 & 5 & 12 & 2 & 10 & 36 & 0.30 & 0.50 & 272 & 66 & 2,338 & 284 & 2,182 \\
miniDAFJS20 & 5 & 13 & 2 & 11 & 40 & 0.54 & 0.52 & 350 & 72 & 3,447 & 363 & 3,274 \\
miniDAFJS21 & 5 & 19 & 3 & 20 & 60 & 0.52 & 0.54 & 730 & 104 & 9,298 & 749 & 9,039 \\
miniDAFJS22 & 5 & 18 & 3 & 18 & 57 & 0.22 & 0.54 & 665 & 99 & 8,223 & 683 & 7,977 \\
miniDAFJS23 & 5 & 21 & 3 & 18 & 66 & 0.53 & 0.54 & 874 & 114 & 11,934 & 895 & 11,649 \\
miniDAFJS24 & 5 & 18 & 3 & 17 & 56 & 0.36 & 0.53 & 642 & 98 & 7,809 & 660 & 7,567 \\
miniDAFJS25 & 5 & 17 & 3 & 17 & 51 & 0.12 & 0.50 & 535 & 91 & 6,024 & 552 & 5,803 \\
miniDAFJS26 & 5 & 17 & 3 & 15 & 53 & 0.36 & 0.53 & 573 & 93 & 6,560 & 590 & 6,331 \\
miniDAFJS27 & 5 & 19 & 3 & 16 & 57 & 0.46 & 0.50 & 721 & 101 & 10,155 & 740 & 9,908 \\
miniDAFJS28 & 5 & 16 & 3 & 15 & 48 & 0.54 & 0.50 & 468 & 86 & 4,871 & 484 & 4,663 \\
miniDAFJS29 & 5 & 14 & 3 & 15 & 47 & 0.33 & 0.59 & 451 & 81 & 4,634 & 465 & 4,432 \\
miniDAFJS30 & 5 & 22 & 3 & 21 & 61 & 0.32 & 0.44 & 791 & 111 & 11,020 & 813 & 10,754 \\
\hline\\
\hline
miniYFJS01 & 7 & 16 & 4 & 12 & 54 & 0.25 & 0.40 & 420 & 94 & 3,554 & 436 & 3,322 \\
miniYFJS02 & 7 & 16 & 4 & 12 & 44 & 0.33 & 0.29 & 294 & 84 & 2,298 & 310 & 2,106 \\
miniYFJS03 & 7 & 16 & 4 & 12 & 45 & 0.08 & 0.30 & 333 & 85 & 2,915 & 349 & 2,719 \\
miniYFJS04 & 7 & 16 & 4 & 12 & 53 & 0.25 & 0.39 & 439 & 93 & 4,083 & 455 & 3,855 \\
miniYFJS05 & 7 & 16 & 4 & 12 & 55 & 0.50 & 0.41 & 457 & 95 & 4,249 & 473 & 4,013 \\
miniYFJS06 & 7 & 16 & 4 & 12 & 51 & 0.42 & 0.36 & 407 & 91 & 3,773 & 423 & 3,553 \\
miniYFJS07 & 7 & 16 & 4 & 12 & 48 & 0.17 & 0.33 & 350 & 88 & 2,912 & 366 & 2,704 \\
miniYFJS08 & 7 & 16 & 4 & 12 & 49 & 0.50 & 0.34 & 363 & 89 & 3,037 & 379 & 2,825 \\
miniYFJS09 & 7 & 16 & 4 & 12 & 51 & 0.33 & 0.36 & 399 & 91 & 3,569 & 415 & 3,349 \\
miniYFJS10 & 7 & 16 & 4 & 12 & 59 & 0.17 & 0.45 & 509 & 99 & 4,767 & 525 & 4,515 \\
miniYFJS11 & 7 & 20 & 5 & 15 & 56 & 0.13 & 0.30 & 464 & 104 & 4,241 & 484 & 3,997 \\
miniYFJS12 & 7 & 20 & 5 & 15 & 68 & 0.13 & 0.40 & 716 & 116 & 8,363 & 736 & 8,071 \\
miniYFJS13 & 7 & 20 & 5 & 15 & 69 & 0.40 & 0.41 & 723 & 117 & 8,272 & 743 & 7,976 \\
miniYFJS14 & 7 & 20 & 5 & 15 & 59 & 0.60 & 0.33 & 509 & 107 & 4,772 & 529 & 4,516 \\
miniYFJS15 & 7 & 20 & 5 & 15 & 53 & 0.33 & 0.28 & 429 & 101 & 3,944 & 449 & 3,712 \\
miniYFJS16 & 7 & 20 & 5 & 15 & 63 & 0.40 & 0.36 & 617 & 111 & 6,838 & 637 & 6,566 \\
miniYFJS17 & 7 & 20 & 5 & 15 & 57 & 0.27 & 0.31 & 485 & 105 & 4,576 & 505 & 4,328 \\
miniYFJS18 & 7 & 20 & 5 & 15 & 51 & 0.07 & 0.26 & 395 & 99 & 3,514 & 415 & 3,290 \\
miniYFJS19 & 7 & 20 & 5 & 15 & 58 & 0.47 & 0.32 & 512 & 106 & 5,031 & 532 & 4,779 \\
miniYFJS20 & 7 & 20 & 5 & 15 & 62 & 0.53 & 0.35 & 638 & 110 & 7,295 & 658 & 7,027 \\
miniYFJS21 & 7 & 24 & 6 & 18 & 72 & 0.28 & 0.33 & 756 & 128 & 8,466 & 780 & 8,154 \\
miniYFJS22 & 7 & 24 & 6 & 18 & 82 & 0.28 & 0.40 & 996 & 138 & 12,962 & 1,020 & 12,610 \\
miniYFJS23 & 7 & 24 & 6 & 18 & 78 & 0.17 & 0.38 & 908 & 134 & 11,280 & 932 & 10,944 \\
miniYFJS24 & 7 & 24 & 6 & 18 & 62 & 0.28 & 0.26 & 588 & 118 & 6,136 & 612 & 5,864 \\
miniYFJS25 & 7 & 24 & 6 & 18 & 76 & 0.17 & 0.36 & 836 & 132 & 9,686 & 860 & 9,358 \\
miniYFJS26 & 7 & 24 & 6 & 18 & 67 & 0.33 & 0.30 & 693 & 123 & 7,841 & 717 & 7,549 \\
miniYFJS27 & 7 & 24 & 6 & 18 & 81 & 0.39 & 0.40 & 951 & 137 & 11,709 & 975 & 11,361 \\
miniYFJS28 & 7 & 24 & 6 & 18 & 67 & 0.22 & 0.30 & 661 & 123 & 7,001 & 685 & 6,709 \\
miniYFJS29 & 7 & 24 & 6 & 18 & 80 & 0.39 & 0.39 & 934 & 136 & 11,536 & 958 & 11,192 \\
miniYFJS30 & 7 & 24 & 6 & 18 & 72 & 0.44 & 0.33 & 800 & 128 & 9,756 & 824 & 9,444 \\
\hline
\end{tabular}
}
\caption{Main features of the proposed sixty small-sized instances.}
\label{tab1}
\end{table}

\begin{table}[ht!]
\centering
\resizebox{0.8\textwidth}{!}{
\rowcolors{3}{}{lightgray}
\begin{tabular}{ccccccccrrrrr} 
\hline
\multicolumn{8}{c}{Main instance characteristics} &
\multicolumn{3}{c}{MILP formulation} &
\multicolumn{2}{c}{CP Optimizer formulation} \\ 
\hline
Instance & $|\mathcal{F}|$ & $|\mathcal{O}|$ & $n$ & $|\widehat A|$ & $\sum_{i \in \mathcal{O}} \mathcal{F}_i$ & $\omega_1$ & $\omega_2$ &
\begin{tabular}[c]{@{}r@{}}\#binary\\ variables\end{tabular} &
\begin{tabular}[c]{@{}r@{}}\#continous\\ variables\end{tabular} & \#constraints &
\begin{tabular}[c]{@{}r@{}}\#interval\\ variables\end{tabular} & \#constraints \\ 
\hline
DAFJS01 & 5 & 26 & 4 & 26 & 82 & 0.54 & 0.54 & 1,358 & 140 & 23,102 & 1,384 & 22,748 \\
DAFJS02 & 5 & 25 & 4 & 23 & 79 & 0.45 & 0.54 & 1,273 & 135 & 21,348 & 1,298 & 21,007 \\
DAFJS03 & 10 & 55 & 4 & 52 & 279 & 0.32 & 0.45 & 7,849 & 400 & 223,911 & 7,904 & 222,740 \\
DAFJS04 & 10 & 43 & 4 & 40 & 220 & 0.25 & 0.46 & 4,960 & 317 & 115,208 & 5,003 & 114,285 \\
DAFJS05 & 5 & 39 & 6 & 34 & 104 & 0.35 & 0.42 & 2,242 & 188 & 50,228 & 2,281 & 49,773 \\
DAFJS06 & 5 & 44 & 6 & 41 & 136 & 0.38 & 0.52 & 3,724 & 230 & 103,253 & 3,768 & 102,665 \\
DAFJS07 & 10 & 85 & 6 & 82 & 431 & 0.30 & 0.45 & 18,695 & 612 & 817,681 & 18,780 & 815,872 \\
DAFJS08 & 10 & 85 & 6 & 82 & 403 & 0.31 & 0.42 & 16,357 & 584 & 670,739 & 16,442 & 669,042 \\
DAFJS09 & 5 & 45 & 8 & 42 & 135 & 0.40 & 0.50 & 3,755 & 231 & 107,667 & 3,800 & 107,082 \\
DAFJS10 & 5 & 58 & 8 & 52 & 168 & 0.40 & 0.47 & 5,764 & 290 & 201,720 & 5,822 & 200,990 \\
DAFJS11 & 10 & 113 & 8 & 108 & 534 & 0.40 & 0.41 & 28,648 & 771 & 1,546,720 & 28,761 & 1,544,471 \\
DAFJS12 & 10 & 117 & 8 & 114 & 603 & 0.49 & 0.46 & 36,513 & 848 & 2,223,141 & 36,630 & 2,220,612 \\
DAFJS13 & 5 & 62 & 10 & 55 & 193 & 0.41 & 0.53 & 7,511 & 323 & 295,582 & 7,573 & 294,748 \\
DAFJS14 & 5 & 69 & 10 & 62 & 206 & 0.37 & 0.50 & 8,578 & 350 & 361,926 & 8,647 & 361,033 \\
DAFJS15 & 10 & 120 & 10 & 117 & 595 & 0.32 & 0.44 & 35,811 & 846 & 2,184,302 & 35,931 & 2,181,802 \\
DAFJS16 & 10 & 120 & 10 & 114 & 602 & 0.33 & 0.45 & 36,344 & 853 & 2,203,066 & 36,464 & 2,200,538 \\
DAFJS17 & 5 & 82 & 12 & 77 & 246 & 0.43 & 0.50 & 12,244 & 416 & 617,041 & 12,326 & 615,975 \\
DAFJS18 & 5 & 74 & 12 & 64 & 231 & 0.41 & 0.53 & 10,785 & 385 & 509,435 & 10,859 & 508,437 \\
DAFJS19 & 7 & 70 & 8 & 66 & 283 & 0.34 & 0.51 & 11,507 & 431 & 471,949 & 11,577 & 470,747 \\
DAFJS20 & 7 & 92 & 10 & 87 & 361 & 0.36 & 0.49 & 18,709 & 553 & 976,026 & 18,801 & 974,490 \\
DAFJS21 & 7 & 107 & 12 & 102 & 425 & 0.38 & 0.50 & 25,853 & 647 & 1,577,783 & 25,960 & 1,575,976 \\
DAFJS22 & 7 & 116 & 12 & 109 & 450 & 0.39 & 0.48 & 29,296 & 690 & 1,932,281 & 29,412 & 1,930,365 \\
DAFJS23 & 9 & 76 & 8 & 71 & 367 & 0.31 & 0.48 & 15,103 & 529 & 628,938 & 15,179 & 627,394 \\
DAFJS24 & 9 & 92 & 8 & 87 & 463 & 0.31 & 0.50 & 23,893 & 657 & 1,238,922 & 23,985 & 1,236,978 \\
DAFJS25 & 9 & 123 & 10 & 119 & 619 & 0.31 & 0.50 & 42,753 & 875 & 2,967,430 & 42,876 & 2,964,831 \\
DAFJS26 & 9 & 119 & 10 & 116 & 606 & 0.34 & 0.51 & 41,026 & 854 & 2,794,800 & 41,145 & 2,792,257 \\
DAFJS27 & 9 & 127 & 12 & 118 & 625 & 0.27 & 0.49 & 43,461 & 889 & 3,029,045 & 43,588 & 3,026,418 \\
DAFJS28 & 10 & 91 & 8 & 89 & 457 & 0.32 & 0.45 & 21,065 & 650 & 980,790 & 21,156 & 978,871 \\
DAFJS29 & 10 & 95 & 8 & 94 & 468 & 0.34 & 0.44 & 22,450 & 669 & 1,109,378 & 22,545 & 1,107,411 \\
DAFJS30 & 10 & 98 & 10 & 94 & 509 & 0.20 & 0.47 & 26,059 & 716 & 1,344,045 & 26,157 & 1,341,911 \\
\hline\\
\hline
YFJS01 & 7 & 40 & 4 & 36 & 104 & 0.10 & 0.27 & 1,824 & 192 & 38,286 & 1,864 & 37,830 \\
YFJS02 & 7 & 40 & 4 & 36 & 104 & 0.17 & 0.27 & 1,568 & 192 & 24,498 & 1,608 & 24,042 \\
YFJS03 & 7 & 24 & 6 & 18 & 63 & 0.28 & 0.27 & 611 & 119 & 6,561 & 635 & 6,285 \\
YFJS04 & 7 & 28 & 7 & 21 & 71 & 0.19 & 0.26 & 813 & 135 & 10,260 & 841 & 9,948 \\
YFJS05 & 7 & 32 & 8 & 24 & 81 & 0.33 & 0.26 & 1,003 & 153 & 13,417 & 1,035 & 13,061 \\
YFJS06 & 7 & 36 & 9 & 27 & 95 & 0.19 & 0.27 & 1,365 & 175 & 21,214 & 1,401 & 20,798 \\
YFJS07 & 7 & 36 & 9 & 27 & 93 & 0.26 & 0.26 & 1,279 & 173 & 18,588 & 1,315 & 18,180 \\
YFJS08 & 12 & 36 & 9 & 27 & 100 & 0.26 & 0.16 & 888 & 185 & 8,879 & 924 & 8,443 \\
YFJS09 & 12 & 36 & 9 & 27 & 219 & 0.22 & 0.46 & 4,079 & 304 & 78,474 & 4,115 & 77,562 \\
YFJS10 & 12 & 40 & 10 & 30 & 113 & 0.17 & 0.17 & 1,169 & 206 & 13,593 & 1,209 & 13,101 \\
YFJS11 & 10 & 50 & 10 & 40 & 134 & 0.22 & 0.19 & 1,860 & 245 & 27,378 & 1,910 & 26,792 \\
YFJS12 & 10 & 50 & 10 & 40 & 133 & 0.12 & 0.18 & 1,915 & 244 & 30,085 & 1,965 & 29,503 \\
YFJS13 & 10 & 50 & 10 & 40 & 137 & 0.15 & 0.19 & 1,895 & 248 & 27,177 & 1,945 & 26,579 \\
YFJS14 & 26 & 221 & 13 & 208 & 641 & 0.24 & 0.08 & 16,603 & 1,110 & 454,179 & 16,824 & 451,394 \\
YFJS15 & 26 & 221 & 13 & 208 & 648 & 0.23 & 0.08 & 16,620 & 1,117 & 441,752 & 16,841 & 438,939 \\
YFJS16 & 26 & 221 & 13 & 208 & 633 & 0.13 & 0.07 & 16,037 & 1,102 & 424,253 & 16,258 & 421,500 \\
YFJS17 & 26 & 289 & 17 & 272 & 1328 & 0.15 & 0.14 & 68,502 & 1,933 & 3,572,342 & 68,791 & 3,566,741 \\
YFJS18 & 26 & 289 & 17 & 272 & 1362 & 0.15 & 0.15 & 72,354 & 1,967 & 3,901,378 & 72,643 & 3,895,641 \\
YFJS19 & 26 & 289 & 17 & 272 & 1347 & 0.20 & 0.15 & 70,527 & 1,952 & 3,737,737 & 70,816 & 3,732,060 \\
YFJS20 & 26 & 289 & 17 & 272 & 1343 & 0.12 & 0.15 & 70,371 & 1,948 & 3,742,757 & 70,660 & 3,737,096 \\
\hline
\end{tabular}
}
\caption{Main features of the fifty large-sized instances from~\cite{birgin2014milp}.}
\label{tab2}
\end{table}

\section{Numerical experiments} \label{experiments}

In this section we present numerical experiments. First, we wish to
evaluate the two introduced constructive heuristics. Second, we wish
to assess the correctness of the MILP and CP models and attempt to
infer which of the two, or rather which of the exact commercial
solvers applied to each of them, is more effective in finding proven
optimal solutions. Third, we wish to determine the usefulness of
providing a feasible solution to the exact solvers. It should be noted
that all efforts are to build a set of test instances with proven
optimal solutions. The models and constructive heuristics presented in
this paper are intended to contribute in that respect and are not
intended to construct a solution method per se, for a known difficult
problem. In all cases we considered the~110 instances introduced in
Section~\ref{instances} with the learning rate $\alpha \in \{0.1, 0.2,
0.3 \}$ for a total of~330 instances.

The experiments were carried out in an Intel i9-12900K (12th Gen)
processor operating at 5.200GHz and 128 GB of RAM. The constructive
heuristics were implemented in C++ programming language. Models were
solved using IBM ILOG CPLEX Optimization Studio version 22.1, using
default parameters, with concert library and C++. The code was
compiled using g++ 10.2.1. Benchmark instances and code are available
at \url{https://github.com/kennedy94/FJS}.

\subsection{Experiments with the constructive heuristics}

In this subsection, we evaluate the two constructive
heuristics. Tables~\ref{tab3} and~\ref{tab4} show the results. For
each instance, the lowest makespan, among the solutions found by the
two constructive heuristics, is shown in bold. In all instances, the
constructive heuristics take less than 0.001 seconds of CPU time to
build a solution. Considering the small-sized instances in
Table~\ref{tab3}, depending on the instance, there may be a
significant difference between the solutions found by one and the
other constructive heuristic. However, on average, the two heuristics
behave basically indistinguishably. For the large-sized instances in
Table~\ref{tab4} the comparison is a bit different: there is a clear
advantage of the EST constructive heuristic in the DA-type instances,
while on the other hand there is a clear advantage of the ECT
constructive heuristic in the Y-type instances. It is worth noting
that in the small-sized instances there is no clear differentiation
between the sequencing flexibility sparsity measure~$\omega_1$ of the
DA-type and the Y-type instances (see Table~\ref{tab1}). On the other
hand, in the large-sized instances, Y-type instances have a value
of~$\omega_1$ clearly lower than the value of~$\omega_1$ of DA-type
instances. Summarizing, as mentioned at the end of
Section~\ref{heuristics}, the greedy strategy of ECT of choosing the
operation/machine pair that terminates first seems to compensate in
situations where, because there is already little sequencing
flexibility, the greedy choice does not cause a large decrease of the
search space.

\begin{table}[ht!]
\centering
\resizebox{!}{0.5\textheight}{
\rowcolors{3}{}{lightgray}
\begin{tabular}{cccccccccc}
\hline
& & \multicolumn{2}{c}{$\alpha = 0.1$}
& & \multicolumn{2}{c}{$\alpha = 0.2$}
& & \multicolumn{2}{c}{$\alpha = 0.3$} \\
\cline{3-4} \cline{6-7} \cline{9-10} 
\multirow{-2}{*}{instance} & & EST & ECT & & EST & ECT & & EST & ECT \\
\cline{1-1} \cline{3-4} \cline{6-7} \cline{9-10} 
miniDAFJS01 &  & \textbf{23,264} & 23,948 &  & \textbf{22,156} & 23,199 &  & \textbf{20,757} & 21,535 \\
miniDAFJS02 &  & \textbf{23,242} & 23,860 &  & \textbf{22,161} & 22,775 &  & \textbf{21,152} & 21,461 \\
miniDAFJS03 &  & \textbf{18,363} & \textbf{18,363} &  & \textbf{17,972} & \textbf{17,972} &  & \textbf{17,621} & \textbf{17,621} \\
miniDAFJS04 &  & \textbf{21,690} & \textbf{21,690} &  & 21,233 & \textbf{21,121} &  & 20,824 & \textbf{20,391} \\
miniDAFJS05 &  & 24,279 & \textbf{22,598} &  & 21,815 & \textbf{20,826} &  & 20,623 & \textbf{19,253} \\
miniDAFJS06 &  & 23,726 & \textbf{23,370} &  & 23,781 & \textbf{22,193} &  & \textbf{21,148} & 21,480 \\
miniDAFJS07 &  & \textbf{28,644} & 33,413 &  & 27,551 & \textbf{25,088} &  & 26,600 & \textbf{24,767} \\
miniDAFJS08 &  & \textbf{19,878} & 25,100 &  & \textbf{18,857} & 22,921 &  & \textbf{17,927} & 20,943 \\
miniDAFJS09 &  & \textbf{25,425} & 27,806 &  & \textbf{23,830} & 25,400 &  & \textbf{22,398} & 23,295 \\
miniDAFJS10 &  & 22,847 & \textbf{21,563} &  & 21,730 & \textbf{21,021} &  & 20,349 & \textbf{20,176} \\
miniDAFJS11 &  & 35,166 & \textbf{33,911} &  & 32,214 & \textbf{31,269} &  & 29,739 & \textbf{28,718} \\
miniDAFJS12 &  & \textbf{20,342} & 21,435 &  & \textbf{19,624} & 20,370 &  & \textbf{19,094} & 19,394 \\
miniDAFJS13 &  & \textbf{18,313} & 18,681 &  & \textbf{17,143} & 17,181 &  & 16,077 & \textbf{15,425} \\
miniDAFJS14 &  & \textbf{26,503} & 27,525 &  & \textbf{24,958} & 25,919 &  & \textbf{23,552} & 24,467 \\
miniDAFJS15 &  & 24,961 & \textbf{22,472} &  & 23,710 & \textbf{21,969} &  & 22,545 & \textbf{20,524} \\
miniDAFJS16 &  & 25,845 & \textbf{25,691} &  & 24,699 & \textbf{24,593} &  & 26,750 & \textbf{23,797} \\
miniDAFJS17 &  & \textbf{21,070} & 21,144 &  & \textbf{20,519} & \textbf{20,519} &  & \textbf{19,941} & \textbf{19,941} \\
miniDAFJS18 &  & 19,716 & \textbf{19,308} &  & \textbf{18,829} & 19,133 &  & \textbf{17,784} & 18,395 \\
miniDAFJS19 &  & \textbf{21,293} & 21,329 &  & \textbf{20,107} & 20,340 &  & \textbf{19,030} & 19,426 \\
miniDAFJS20 &  & \textbf{23,570} & 26,421 &  & 21,759 & \textbf{21,286} &  & 20,140 & \textbf{19,587} \\
miniDAFJS21 &  & 25,764 & \textbf{25,020} &  & 24,040 & \textbf{22,689} &  & 23,171 & \textbf{21,299} \\
miniDAFJS22 &  & \textbf{28,122} & 31,892 &  & \textbf{26,289} & 29,188 &  & \textbf{24,671} & 26,829 \\
miniDAFJS23 &  & \textbf{27,566} & 33,389 &  & \textbf{25,592} & 30,242 &  & \textbf{23,847} & 28,000 \\
miniDAFJS24 &  & \textbf{26,932} & 27,581 &  & 29,407 & \textbf{24,500} &  & 27,040 & \textbf{22,524} \\
miniDAFJS25 &  & \textbf{23,370} & 25,879 &  & \textbf{22,613} & 24,209 &  & \textbf{21,923} & 22,976 \\
miniDAFJS26 &  & 26,408 & \textbf{25,974} &  & 24,476 & \textbf{23,630} &  & 22,771 & \textbf{21,613} \\
miniDAFJS27 &  & \textbf{28,084} & 31,060 &  & \textbf{25,429} & 26,652 &  & \textbf{23,089} & 23,595 \\
miniDAFJS28 &  & \textbf{28,836} & 30,947 &  & \textbf{26,587} & 27,477 &  & 26,831 & \textbf{24,429} \\
miniDAFJS29 &  & 27,569 & \textbf{21,151} &  & 25,923 & \textbf{19,789} &  & 24,444 & \textbf{18,560} \\
miniDAFJS30 &  & 27,870 & \textbf{27,249} &  & 25,318 & \textbf{24,491} &  & \textbf{21,421} & 22,071 \\ \hline
wins &  & 19 & 13 &  & 16 & 16 &  & 16 & 16 \\
mean &  & 24,621.93 & 25,325.67 &  & 23,344.07 & 23,265.40 &  & 22,108.63 & 21,749.73 \\ \hline
\rowcolor{white} \\
\hline
\rowcolor{white}
& & \multicolumn{2}{c}{$\alpha = 0.1$}
& & \multicolumn{2}{c}{$\alpha = 0.2$}
& & \multicolumn{2}{c}{$\alpha = 0.3$} \\
\cline{3-4} \cline{6-7} \cline{9-10} 
\multirow{-2}{*}{instance} & & EST & ECT & & EST & ECT & & EST & ECT \\
  \cline{1-1} \cline{3-4} \cline{6-7} \cline{9-10} 
miniYFJS01 &  & 40,456 & \textbf{35,243} &  & 38,132 & \textbf{34,443} &  & 36,008 & \textbf{33,697} \\
miniYFJS02 &  & 34,009 & \textbf{28,688} &  & 31,940 & \textbf{27,557} &  & 27,987 & \textbf{25,969} \\
miniYFJS03 &  & \textbf{62,704} & 69,005 &  & \textbf{57,950} & 60,302 &  & 53,679 & \textbf{52,865} \\
miniYFJS04 &  & 32,133 & \textbf{25,394} &  & 30,048 & \textbf{24,669} &  & 24,428 & \textbf{24,017} \\
miniYFJS05 &  & 29,009 & \textbf{27,650} &  & 26,828 & \textbf{26,135} &  & 27,047 & \textbf{24,743} \\
miniYFJS06 &  & 49,001 & \textbf{30,952} &  & 43,387 & \textbf{29,080} &  & 40,871 & \textbf{27,366} \\
miniYFJS07 &  & \textbf{51,782} & 57,652 &  & \textbf{47,743} & 53,702 &  & \textbf{44,128} & 50,110 \\
miniYFJS08 &  & \textbf{35,336} & 51,013 &  & \textbf{33,178} & 40,884 &  & \textbf{31,133} & 40,567 \\
miniYFJS09 &  & \textbf{38,626} & 39,463 &  & 40,411 & \textbf{37,593} &  & 38,480 & \textbf{34,357} \\
miniYFJS10 &  & \textbf{32,975} & 34,363 &  & \textbf{30,651} & 31,742 &  & \textbf{28,513} & 29,397 \\
miniYFJS11 &  & 51,600 & \textbf{51,212} &  & 48,172 & \textbf{47,079} &  & 45,079 & \textbf{43,254} \\
miniYFJS12 &  & \textbf{38,889} & 39,024 &  & \textbf{36,036} & 36,691 &  & \textbf{33,495} & 35,223 \\
miniYFJS13 &  & 38,448 & \textbf{35,177} &  & 34,046 & \textbf{31,720} &  & \textbf{26,536} & 28,751 \\
miniYFJS14 &  & \textbf{41,560} & 42,764 &  & \textbf{39,254} & 39,862 &  & \textbf{37,157} & 37,177 \\
miniYFJS15 &  & \textbf{57,017} & 58,693 &  & \textbf{51,504} & \textbf{51,504} &  & \textbf{46,583} & \textbf{46,583} \\
miniYFJS16 &  & 39,097 & \textbf{35,831} &  & 36,810 & \textbf{33,285} &  & 34,799 & \textbf{30,576} \\
miniYFJS17 &  & \textbf{53,052} & 62,576 &  & \textbf{48,745} & 57,600 &  & \textbf{44,916} & 50,137 \\
miniYFJS18 &  & \textbf{37,217} & 39,284 &  & 34,053 & \textbf{31,911} &  & 31,154 & \textbf{29,764} \\
miniYFJS19 &  & 46,640 & \textbf{46,478} &  & 41,729 & \textbf{41,451} &  & 42,825 & \textbf{37,034} \\
miniYFJS20 &  & \textbf{42,408} & 42,465 &  & 39,852 & \textbf{39,600} &  & 37,509 & \textbf{37,063} \\
miniYFJS21 &  & 50,076 & \textbf{44,016} &  & 45,541 & \textbf{40,765} &  & 40,884 & \textbf{35,905} \\
miniYFJS22 &  & \textbf{43,722} & 47,389 &  & 40,428 & \textbf{32,485} &  & 31,887 & \textbf{30,212} \\
miniYFJS23 &  & \textbf{53,046} & 59,584 &  & \textbf{49,209} & 52,628 &  & \textbf{45,823} & 46,589 \\
miniYFJS24 &  & \textbf{42,268} & 45,974 &  & \textbf{39,883} & 41,759 &  & \textbf{37,651} & 37,998 \\
miniYFJS25 &  & \textbf{46,936} & 68,698 &  & \textbf{40,933} & 61,260 &  & \textbf{35,966} & 54,247 \\
miniYFJS26 &  & 55,031 & \textbf{55,022} &  & 49,428 & \textbf{49,412} &  & \textbf{53,604} & 54,006 \\
miniYFJS27 &  & 43,319 & \textbf{41,060} &  & 39,256 & \textbf{36,561} &  & 35,655 & \textbf{33,799} \\
miniYFJS28 &  & 46,385 & \textbf{43,232} &  & 42,230 & \textbf{40,325} &  & 39,363 & \textbf{37,547} \\
miniYFJS29 &  & \textbf{48,806} & 55,473 &  & 43,757 & \textbf{41,834} &  & 39,279 & \textbf{39,146} \\
miniYFJS30 &  & \textbf{49,458} & 52,043 &  & \textbf{45,051} & 46,003 &  & 41,116 & \textbf{40,752} \\ \hline
wins &  & 17 & 13 &  & 12 & 19 &  & 12 & 19 \\
mean &  & 44,366.87 & 45,513.93 &  & 40,872.83 & 40,661.40 &  & 37,785.17 & 37,628.37 \\ \hline
\end{tabular}}
\caption{Makespan values for the small-sized set of instance solved by constructive heuristics.}
\label{tab3}
\end{table}

\begin{table}[ht!]
\centering
\resizebox{!}{0.5\textheight}{
\rowcolors{3}{}{lightgray}
\begin{tabular}{cccccccccc}
\hline
& & \multicolumn{2}{c}{$\alpha = 0.1$}
& & \multicolumn{2}{c}{$\alpha = 0.2$}
& & \multicolumn{2}{c}{$\alpha = 0.3$} \\
\cline{3-4} \cline{6-7} \cline{9-10} 
\multirow{-2}{*}{instance} & & EST & ECT & & EST & ECT & & EST & ECT \\
\cline{1-1} \cline{3-4} \cline{6-7} \cline{9-10}
DAFJS01 & & \textbf{29,769} & 41,358 & & \textbf{28,920} & 31,657 & & \textbf{22,616} & 28,222 \\
DAFJS02 & & \textbf{32,467} & 33,155 & & 29,089 & \textbf{28,154} & & \textbf{26,183} & 26,949 \\
DAFJS03 & & \textbf{53,688} & 53,834 & & \textbf{48,555} & 51,717 & & \textbf{43,964} & 45,501 \\
DAFJS04 & & \textbf{54,082} & 54,150 & & \textbf{48,461} & 48,503 & & \textbf{43,281} & 44,094 \\
DAFJS05 & & 52,651 & \textbf{45,862} & & \textbf{44,790} & 44,919 & & \textbf{38,237} & 38,457 \\
DAFJS06 & & \textbf{51,925} & 56,012 & & \textbf{44,228} & 46,403 & & \textbf{38,007} & 41,229 \\
DAFJS07 & & \textbf{57,193} & 59,491 & & 51,019 & \textbf{50,802} & & \textbf{41,285} & 44,518 \\
DAFJS08 & & \textbf{62,159} & 67,107 & & 53,998 & \textbf{53,439} & & \textbf{46,323} & 47,173 \\
DAFJS09 & & \textbf{48,680} & 60,565 & & 48,833 & \textbf{42,537} & & \textbf{41,418} & 42,518 \\
DAFJS10 & & \textbf{58,695} & 60,374 & & \textbf{46,257} & 49,140 & & \textbf{38,771} & 41,125 \\
DAFJS11 & & \textbf{67,594} & 92,639 & & \textbf{56,847} & 75,992 & & \textbf{47,805} & 56,547 \\
DAFJS12 & & 70,287 & \textbf{68,322} & & 58,100 & \textbf{55,575} & & \textbf{45,735} & 49,991 \\
DAFJS13 & & 63,386 & \textbf{61,794} & & \textbf{51,317} & 53,486 & & \textbf{42,319} & 44,724 \\
DAFJS14 & & 83,362 & \textbf{76,862} & & \textbf{62,724} & 63,887 & & \textbf{49,459} & 50,737 \\
DAFJS15 & & 78,413 & \textbf{72,288} & & \textbf{54,353} & 69,279 & & \textbf{49,206} & 54,179 \\
DAFJS16 & & 78,289 & \textbf{76,204} & & \textbf{65,550} & 72,885 & & 57,213 & \textbf{56,727} \\
DAFJS17 & & \textbf{73,219} & 84,719 & & \textbf{63,177} & 65,307 & & \textbf{54,527} & 54,887 \\
DAFJS18 & & 82,129 & \textbf{78,862} & & \textbf{61,831} & 67,753 & & \textbf{51,823} & 52,325 \\
DAFJS19 & & \textbf{66,412} & 67,169 & & \textbf{58,001} & 63,933 & & 49,746 & \textbf{46,880} \\
DAFJS20 & & \textbf{78,778} & 80,562 & & \textbf{63,588} & 68,034 & & \textbf{51,686} & 54,531 \\
DAFJS21 & & \textbf{78,320} & 83,933 & & \textbf{63,202} & 66,878 & & \textbf{53,010} & 54,906 \\
DAFJS22 & & 77,853 & \textbf{70,892} & & \textbf{56,115} & 62,199 & & \textbf{45,005} & 49,723 \\
DAFJS23 & & \textbf{49,969} & 53,123 & & \textbf{47,616} & 50,639 & & \textbf{39,249} & 41,475 \\
DAFJS24 & & \textbf{57,411} & 62,038 & & \textbf{49,019} & 50,345 & & \textbf{44,851} & 46,119 \\
DAFJS25 & & 89,248 & \textbf{82,055} & & \textbf{66,262} & 67,533 & & \textbf{51,964} & 60,266 \\
DAFJS26 & & \textbf{81,480} & 83,635 & & 75,230 & \textbf{72,226} & & \textbf{57,182} & 60,803 \\
DAFJS27 & & \textbf{81,470} & 88,629 & & \textbf{64,797} & 71,903 & & \textbf{58,138} & 58,145 \\
DAFJS28 & & \textbf{62,568} & 64,560 & & \textbf{52,639} & 53,110 & & \textbf{42,898} & 46,294 \\
DAFJS29 & & 74,841 & \textbf{72,938} & & \textbf{59,792} & 66,900 & & \textbf{51,606} & 58,809 \\
DAFJS30 & & \textbf{61,147} & 70,062 & & \textbf{55,015} & 67,245 & & \textbf{43,859} & 47,605 \\
\hline
wins & & 20 & 10 & & 24 & 6 & & 28 & 2 \\
mean & & 65,249.50 & 67,439.80 & & 54,310.83 & 57,746.00 & & 45,578.87 & 48,181.97 \\
\hline
\rowcolor{white} \\
\hline
\rowcolor{white}
& & \multicolumn{2}{c}{$\alpha = 0.1$}
& & \multicolumn{2}{c}{$\alpha = 0.2$}
& & \multicolumn{2}{c}{$\alpha = 0.3$} \\
\cline{3-4} \cline{6-7} \cline{9-10} 
\rowcolor{white}
\multirow{-2}{*}{instance} & & EST & ECT & & EST & ECT & & EST & ECT \\
\cline{1-1} \cline{3-4} \cline{6-7} \cline{9-10}
YFJS01 & & \textbf{87,203} & 106,117 & & \textbf{84,152} & 92,107 & & \textbf{70,402} & 80,419 \\
YFJS02 & & 87,462 & \textbf{81,579} & & 73,957 & \textbf{66,853} & & 68,075 & \textbf{61,111} \\
YFJS03 & & 42,457 & \textbf{40,197} & & \textbf{35,380} & 37,159 & & \textbf{31,680} & 33,077 \\
YFJS04 & & \textbf{47,467} & 50,724 & & \textbf{43,614} & 44,467 & & 39,249 & \textbf{38,898} \\
YFJS05 & & \textbf{46,138} & 55,871 & & \textbf{40,851} & 49,893 & & \textbf{41,586} & 45,012 \\
YFJS06 & & 53,210 & \textbf{52,487} & & \textbf{48,564} & 54,660 & & \textbf{46,811} & 47,173 \\
YFJS07 & & 63,320 & \textbf{54,457} & & 56,313 & \textbf{51,261} & & \textbf{39,896} & 44,004 \\
YFJS08 & & 51,818 & \textbf{49,626} & & 47,218 & \textbf{44,074} & & \textbf{39,325} & 39,469 \\
YFJS09 & & 38,836 & \textbf{28,354} & & 36,123 & \textbf{26,072} & & 32,120 & \textbf{24,027} \\
YFJS10 & & \textbf{42,583} & 59,808 & & \textbf{40,658} & 52,811 & & \textbf{39,607} & 48,821 \\
YFJS11 & & 65,011 & \textbf{59,356} & & 58,106 & \textbf{51,469} & & 51,943 & \textbf{45,051} \\
YFJS12 & & \textbf{70,830} & 74,978 & & \textbf{59,898} & 62,299 & & 53,275 & \textbf{52,398} \\
YFJS13 & & 53,601 & \textbf{50,805} & & 48,729 & \textbf{45,084} & & 41,303 & \textbf{40,119} \\
YFJS14 & & 151,365 & \textbf{129,428} & & 116,157 & \textbf{109,469} & & 111,313 & \textbf{92,457} \\
YFJS15 & & 152,375 & \textbf{138,196} & & 120,600 & \textbf{112,424} & & 107,006 & \textbf{96,547} \\
YFJS16 & & 144,976 & \textbf{127,055} & & 131,230 & \textbf{106,855} & & 118,515 & \textbf{92,811} \\
YFJS17 & & 133,982 & \textbf{109,112} & & 110,203 & \textbf{85,736} & & 98,045 & \textbf{73,682} \\
YFJS18 & & 154,214 & \textbf{133,703} & & 121,563 & \textbf{99,429} & & 104,338 & \textbf{87,059} \\
YFJS19 & & 133,142 & \textbf{107,055} & & 110,125 & \textbf{89,561} & & 91,287 & \textbf{74,431} \\
YFJS20 & & 137,326 & \textbf{97,868} & & 104,036 & \textbf{91,958} & & 91,229 & \textbf{72,481} \\
\hline
wins & & 5 & 15 & & 7 & 13 & & 7 & 13 \\
mean & & 87,865.80 & 80,338.80 & & 74,373.85 & 68,682.05 && 65,850.25& 59,452.35\\
\hline
\end{tabular}}
\caption{Makespan values for the testbed set of instance solved by constructive heuristics.}
\label{tab4}
\end{table}

\subsection{Solving the proposed models with a commercial solver}

In this section, we apply an exact solver to the~330 instances
considered and evaluate the quality of the solutions found within a
given CPU time limit, depending on whether we provide the solution
found by the constructive heuristics as an initial solution or
not. Since there is no clear winner between the two constructive
heuristics and their execution time is negligible, we run both and
give the best of the two as initial solution to the exact solver. In
the tables and figures, the solvers' runs that receive as input a
solution computed by one of the constructive heuristics are identified
with the expression ``warm start''.

A solution is reported as optimal by the solvers when
\begin{equation}
\text{absolute gap} = \text{best feasible solution} - \text{best lower bound} \leq \epsilon_{\mathrm{abs}}
\end{equation}
or
\begin{equation}
\text{relative gap} = \dfrac{|\text{best feasible solution} - \text{best lower bound}|}{10^{-10} + |\text{best feasible solution}|}
\leq \epsilon_{\mathrm{rel}},
\end{equation}
where, by default, $\epsilon_{\mathrm{abs}} = 10^{-6}$ and
$\epsilon_{\mathrm{rel}} = 10^{-4}$, and ``best feasible solution''
means the smallest value of the makespan related to a feasible
solution generated by the method. Since the optimal makespan of the
instances considered in this work always assumes an integer value, we
choose to use $\epsilon_{\mathrm{abs}} = 1 - 10^{-6}$ and
$\epsilon_{\mathrm{rel}} = 0$. The choice $\epsilon_{\mathrm{rel}} =
0$ avoids premature stops in a solution that may not be optimal.  The
choice $\epsilon_{\mathrm{abs}} = 1 - 10^{-6}$ allows to stop early
when a relative gap smaller than 1 clearly indicates that the optimal
solution has already been found. A CPU time limit of~1 hour was
imposed. All other solvers parameters were used with their default
values.

\begin{figure}[ht!]
\begin{center}
\begin{tabular}{c}
\includegraphics[scale=0.9]{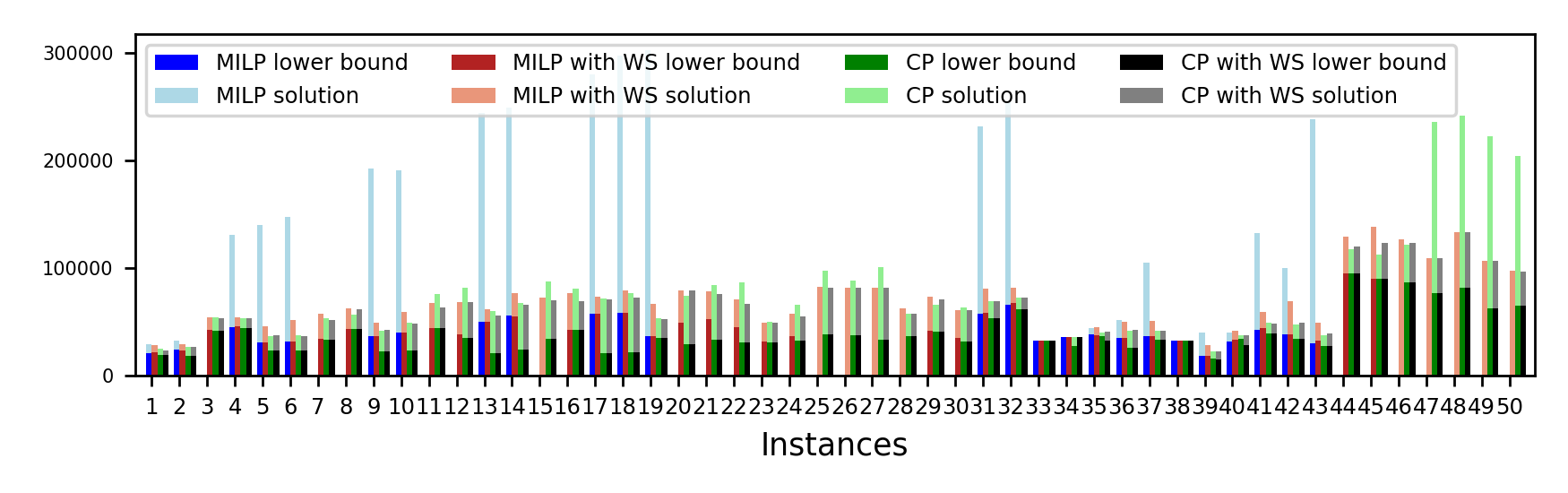}\\
\includegraphics[scale=0.9]{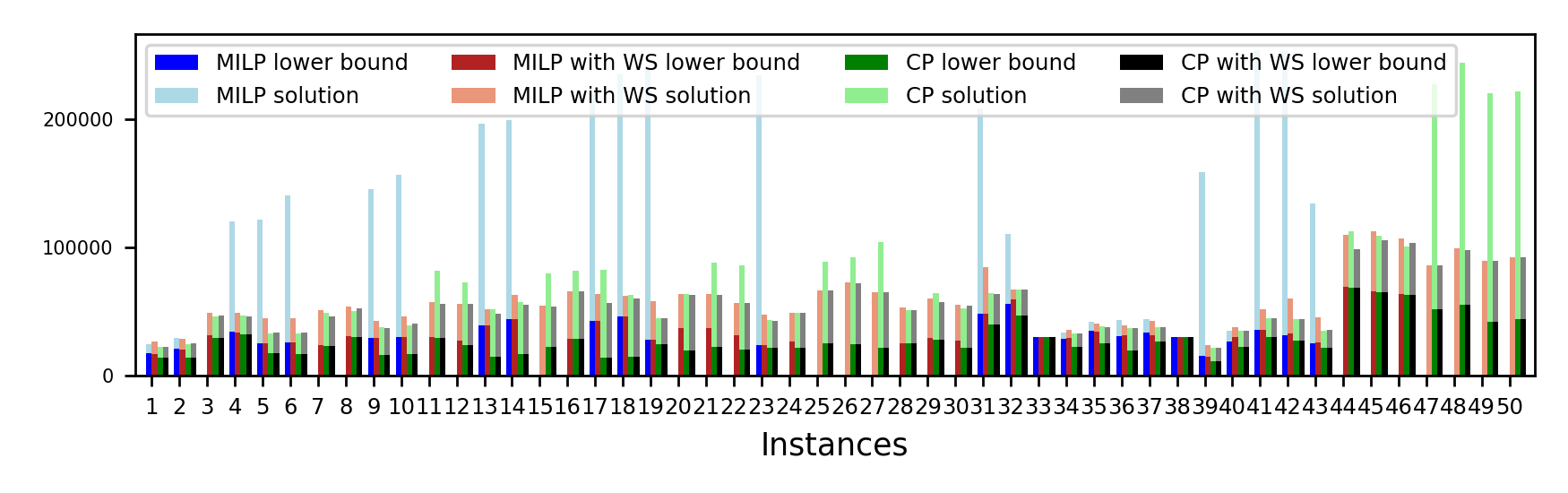}\\
\includegraphics[scale=0.9]{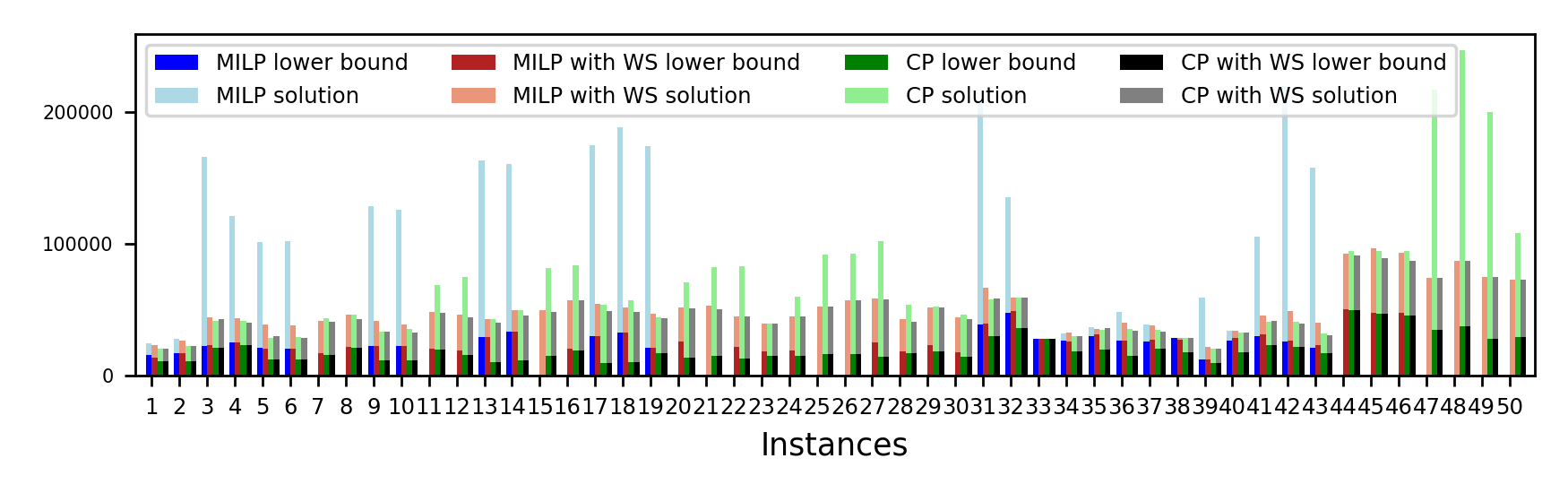}\\
\end{tabular}
\end{center}
\caption{Graphical representation of the information contained in
  Tables~\ref{tab6}--\ref{tab16}. The figures at the top, middle and
  bottom correspond to $\alpha$ equal to $0.1$, $0.2$, and $0.3$,
  respectively. Each figure shows, for each large-sized size
  instance, the solution found by solving the MILP and CP models with
  and without warm start, as well as the lower bounds found in these 4
  cases. The numbers from 1 to 50 on the abscissa axis correspond to
  the 50 large-sized instances in the order presented in all the
  tables. That is, the instances from 1 to 30 correspond to the
  instances DAFJS01, $\dots$, DAFJS30 and the instances from 31 to 50
  correspond to the instances YFJS01, $\dots$, YFJS20.}
\label{fig6}
\end{figure}

Tables~\ref{tab5}, \ref{tab6}, \ref{tab7}, \ref{tab8}, \ref{tab9}, and
\ref{tab10} show the results in detail when the warm start is not
taken into account. Tables~\ref{tab5}, \ref{tab6}, and \ref{tab7}
correspond to the small-sized instances with $\alpha$ equal to $0.1$, $0.2$,
and $0.3$, respectively; while Tables~\ref{tab8}, \ref{tab9}, and
\ref{tab10} show the same thing for the large-sized instances.
Tables~\ref{tab11}--\ref{tab16} show the results when the warm start
is taken into account. The tables show the information related to the
resolution of the MILP and CP models. When a single number appears in
the ``makespan'' column, it means that a provably optimal solution
with that makespan value was found. When instead of a number an
expression of the form $[A,B]\,C\%$ appears, it means that a feasible
solution was found with value~$B$ for the makespan, value~$A$ for a
lower bound on the makespan, and gap~$C$ equal to $100 (B-A)/B$. As
measures of computational effort, for the MILP solver the tables show
the number of iterations, the number of nodes explored in the
branch-and-bound search tree, and the CPU time in seconds. For the CP
solver, the tables show the number of branches and the CPU time in
seconds. If no information is displayed for a particular instance, it
means that the solver was unable to find a feasible solution within
the specified CPU time limit. In Tables~\ref{tab11}--\ref{tab16},
which show results with warm start, the symbol~$\dagger$ next to the
optimal or best value found means that the exact method returned
exactly the same solution computed with the constructive heuristic and
given as initial solution. The information from Tables
~\ref{tab11}--\ref{tab16} is presented graphically in
Figure~\ref{fig6}.

Let's look at the small-sized instances first. Without warm start, the exact methods found provably optimal solutions for 168 instances of the MILP model and 169 instances of the CP model (out of a total of 180 small-sized instances). In the remaining cases, the gaps for the MILP model instances were between 0.17\% and 23.32\%, while for the CP model instances the gaps were between 29.44\% and 47.55\%. It is worth noting that in the few cases without a proven optimal solution, there is a slight advantage in the best solution found for the CP model instances and a slight advantage in the lower bounds found for the MILP model instances. In the instances where a proven solution was found by solving both the MILP model and the CP model, the CP solver was on average 12.09\% faster than the MILP solver. When the constructive heuristics solution is made available to the exact solvers, the number of proven optimal solutions found hardly changes (still the same number for MILP model instances and one less for CP model instances, not necessarily the same as those solved without the warm start). However, for MILP model instances where a proven optimal solution is found both with and without warm start, the use of warm start reduces the solution time by 32.52\% in average. This reduction is 4.44\% for the CP model solver. One way or another, whether using the MILP model or the CP model, with or without warm start, it was possible to find provably optimal solutions in 176 (out of 180) small-sized instances.

The analysis of the~150 large-sized instances is a little different. Without a warm start, the exact methods were able to find proven optimal solutions for only 7 instances of the MILP model and 5 instances of the CP model. In~70 MILP model instances it was possible to find feasible solutions with gaps between 13.94\% and 90.47\%, while feasible solutions with gaps between 9.30\% and 86.10\% were found in 145 instances of the CP model. In 73 instances of the MILP model, not a single feasible solution was found. In the 70 instances in which feasible solutions from both the MILP and CP models were found, the solutions found using the CP model were on average 68.85\% better. It was after observing these results that the idea arose to develop constructive heuristics to test the warm start and consider a set of smaller instances.

When the solution of the constructive heuristics is fed to the exact solver, 6 provably optimal solutions and 144 feasible solutions are obtained with gaps between 5.65\% and 64.36\% for MILP model instances. For the CP model instances, the same number of provably optimal and feasible solutions are found, with gaps between 15.43\% and 81.06\% for the feasible ones. In those instances where a proven optimal solution is found without and with warm start, warm start \textit{increases} the computational cost of solving the MILP and CP models by 9.69\% and 25.99\%, respectively. On the other hand, in the MILP model instances in which a feasible solution was found without the use of warm start, the use of warm start improved the quality of the feasible solution found by 49.46\%. For the CP model instances, this improvement was 11.53\%. In the 144 instances in which feasible solutions from both the MILP and CP models were found, the solutions found using the CP model were on average 6.82\% better. The question remains as to whether the exact methods are able to improve the solution provided by the constructive heuristics or whether the statistics improve only because the solvers return the solution they received as input. In the case of the MILP model instances, the initial solution is improved in 33 problems, while in the CP model instances the initial solution is improved in 134 problems. Without a warm start, in the instances where it is possible to find a provably optimal solution for both the MILP model and the CP model (5 instances), the cost of solving the CP models is 70.81\% lower. In the case where provably optimal solutions are found in both cases using warm start (6 instances), the cost of solving the CP models is 41.86\% lower. In short, it is challenging to find a proven optimal solution for large-sized instances, solving CP model instances costs less, CP models provide better quality feasible solutions when it is not possible to find a provably optimal solution, and solving MILP models provides better quality lower bounds. One way or another, whether using the MILP model or the CP model, with or without warm start, it was possible to find provably optimal solutions in only 7 large-sized instances and feasible solutions in all the remaining 143 large-sized instances.

\begin{table}[ht!]
\centering
\resizebox{0.85\textwidth}{!}{
\rowcolors{7}{}{lightgray}
\begin{tabular}{cccccccccc}
\hline
\multirow{3}{*}{instance}
&  & \multicolumn{4}{c}{IBM ILOG CPLEX}
&  & \multicolumn{3}{c}{IBM ILOG CP Optimizer} \\
\cline{3-6} \cline{8-10} 
&  & \multirow{2}{*}{makespan} & \multicolumn{3}{c}{Effort measurement}
&  & \multirow{2}{*}{makespan} & \multicolumn{2}{c}{Effort measurement} \\
\cline{4-6} \cline{9-10} 
&  &  & \#iterations & \#B\&B Nodes & CPU &  &  & \#branches & CPU \\
\cline{1-1} \cline{3-6} \cline{8-10} 
miniDAFJS01 &  & 22,875 & 278,358 & 2,630 & 20.9 &  & 22,875 & 134,773 & 5.5 \\
miniDAFJS02  &  & 22,708 & 42,653 & 579 & 1.7 &  & 22,708 & 22,910 & 1.1 \\
miniDAFJS03  &  & 18,363 & 10,016 & 153 & 0.9 &  & 18,363 & 14,385 & 0.3 \\
miniDAFJS04  &  & 20,498 & 4,002 & 133 & 0.3 &  & 20,498 & 5,579 & 0.2 \\
miniDAFJS05  &  & 20,593 & 44,305 & 555 & 2.7 &  & 20,593 & 66,626 & 3.0 \\
\hline
\end{tabular}}
\caption{Solutions found and computational cost of solving the
  small-sized instances with learning rate $\alpha = 0.1$ using CPLEX
  and CP Optimizer with no warm start.}
\label{tab5}
\end{table}

\begin{table}[ht!]
\centering
\resizebox{0.85\textwidth}{!}{
\rowcolors{7}{}{lightgray}
\begin{tabular}{cccccccccc}
\hline
\multirow{3}{*}{instance}
&  & \multicolumn{4}{c}{IBM ILOG CPLEX}
&  & \multicolumn{3}{c}{IBM ILOG CP Optimizer} \\
\cline{3-6} \cline{8-10} 
&  & \multirow{2}{*}{makespan} & \multicolumn{3}{c}{Effort measurement}
&  & \multirow{2}{*}{makespan} & \multicolumn{2}{c}{Effort measurement} \\
\cline{4-6} \cline{9-10} 
&  &  & \#iterations & \#B\&B Nodes & CPU &  &  & \#branches & CPU \\
\cline{1-1} \cline{3-6} \cline{8-10} 
miniDAFJS01 &  & 21,327 & 546,297 & 4,285 & 31.5 &  & 21,327 & 953,712 & 51.7 \\
miniDAFJS02  &  & 20,635 & 81,821 & 640 & 3.3 &  & 20,635 & 62,089 & 1.8 \\
miniDAFJS03  &  & 17,972 & 19,340 & 305 & 0.9 &  & 17,972 & 22,808 & 0.5 \\
miniDAFJS04  &  & 19,602 & 7,746 & 226 & 0.3 &  & 19,602 & 13,297 & 0.9 \\
miniDAFJS05  &  & 18,803 & 98,953 & 1,120 & 5.1 &  & 18,803 & 151,224 & 7.9 \\
\hline
\end{tabular}}
\caption{Solutions found and computational cost of solving the
  small-sized instances with learning rate $\alpha = 0.2$ using CPLEX
  and CP Optimizer with no warm start.}
\label{tab6}
\end{table}

\begin{table}[ht!]
\centering
\resizebox{0.85\textwidth}{!}{
\rowcolors{7}{}{lightgray}
\begin{tabular}{cccccccccc}
\hline
\multirow{3}{*}{instance}
&  & \multicolumn{4}{c}{IBM ILOG CPLEX}
&  & \multicolumn{3}{c}{IBM ILOG CP Optimizer} \\
\cline{3-6} \cline{8-10} 
&  & \multirow{2}{*}{makespan} & \multicolumn{3}{c}{Effort measurement}
&  & \multirow{2}{*}{makespan} & \multicolumn{2}{c}{Effort measurement} \\
\cline{4-6} \cline{9-10} 
&  &  & \#iterations & \#B\&B Nodes & CPU &  &  & \#branches & CPU \\
\cline{1-1} \cline{3-6} \cline{8-10} 
miniDAFJS01 &  & 19,443 & 2,072,668 & 16,713 & 114.7 &  & 19,443 & 2,800,465 & 77.8 \\
miniDAFJS02 &  & 18,916 & 53,092 & 696 & 2.5 &  & 18,916 & 59,859 & 2.2 \\
miniDAFJS03 &  & 17,419 & 17,807 & 254 & 1.1 &  & 17,419 & 36,788 & 0.9 \\
miniDAFJS04 &  & 18,800 & 5,444 & 141 & 0.4 &  & 18,800 & 19,010 & 0.7 \\
miniDAFJS05 &  & 17,596 & 149,857 & 1,357 & 26.4 &  & 17,596 & 354,775 & 19.2 \\
\hline
\end{tabular}}
\caption{Solutions found and computational cost of solving the
  small-sized instances with learning rate $\alpha = 0.3$ using CPLEX
  and CP Optimizer with no warm start.}
\label{tab7}
\end{table}

\begin{table}[ht!]
\centering
\resizebox{\textwidth}{!}{
\rowcolors{5}{}{lightgray}
\begin{tabular}{cccccccccc}
\hline
\multirow{3}{*}{instance}
&  & \multicolumn{4}{c}{IBM ILOG CPLEX}
&  & \multicolumn{3}{c}{IBM ILOG CP Optimizer} \\
\cline{3-6} \cline{8-10} 
&  & \multirow{2}{*}{makespan} & \multicolumn{3}{c}{Effort measurement}
&  & \multirow{2}{*}{makespan} & \multicolumn{2}{c}{Effort measurement} \\
\cline{4-6} \cline{9-10} 
&  &  & \#iterations & \#B\&B Nodes & CPU &  &  & \#branches & CPU \\
\cline{1-1} \cline{3-6} \cline{8-10} 
DAFJS01 &  & {[}20,320, 28,850{]} 29.57\% & 2,781,484 & 10,131 & 3,600.0 &  & {[}18,617, 24,424{]} 23.78\% & 8,687,672 & 3,600.0 \\
DAFJS02 &  & {[}23,614, 31,946{]} 26.08\% & 2,850,568 & 10,711 & 3,600.0 &  & {[}18,533, 26,575{]} 30.26\% & 11,632,154 & 3,600.0 \\
DAFJS03 &  &  &  &  & 3,600.0 &  & {[}41,271, 54,098{]} 23.71\% & 801,461 & 3,600.1 \\
DAFJS04 &  & {[}45,286, 131,061{]} 65.45\% & 215,396 & 13 & 3,600.0 &  & {[}44,062, 53,100{]} 17.02\% & 850,903 & 3,600.2 \\
DAFJS05 &  & {[}30,405, 140,075{]} 78.29\% & 291,690 & 510 & 3,600.0 &  & {[}23,378, 36,347{]} 35.68\% & 2,711,057 & 3,600.1 \\
\hline
\end{tabular}}
\caption{Solutions found and computational cost of solving the
  large-sized instances with learning rate $\alpha = 0.1$ using CPLEX
  and CP Optimizer with no warm start.}
\label{tab8}
\end{table}

\begin{table}[ht!]
\centering
\resizebox{\textwidth}{!}{
\rowcolors{5}{}{lightgray}
\begin{tabular}{cccccccccc}
\hline
\multirow{3}{*}{instance}
&  & \multicolumn{4}{c}{IBM ILOG CPLEX}
&  & \multicolumn{3}{c}{IBM ILOG CP Optimizer} \\
\cline{3-6} \cline{8-10} 
&  & \multirow{2}{*}{makespan} & \multicolumn{3}{c}{Effort measurement}
&  & \multirow{2}{*}{makespan} & \multicolumn{2}{c}{Effort measurement} \\
\cline{4-6} \cline{9-10} 
&  &  & \#iterations & \#B\&B Nodes & CPU &  &  & \#branches & CPU \\
\cline{1-1} \cline{3-6} \cline{8-10} 
DAFJS01 &  & {[}17,243, 24,221{]} 28.81\% & 4,060,290 & 15,216 & 3,600.0 &  & {[}14,141, 21,955{]} 35.59\% & 10,060,921 & 3,600.1 \\
DAFJS02 &  & {[}20,460, 29,199{]} 29.93\% & 4,007,884 & 12,017 & 3,600.0 &  & {[}13,966, 24,524{]} 43.05\% & 9,597,884 & 3,600.0 \\
DAFJS03 &  &  &  &  & 3,600.0 &  & {[}29,519, 46,041{]} 35.89\% & 805,535 & 3,600.2 \\
DAFJS04 &  & {[}33,806, 120,313{]} 71.90\% & 227,790 & 80 & 3,600.0 &  & {[}32,038, 46,931{]} 31.73\% & 941,659 & 3,600.4 \\
DAFJS05 &  & {[}25,343, 121,192{]} 79.09\% & 492,382 & 517 & 3,600.0 &  & {[}17,025, 32,708{]} 47.95\% & 2,914,026 & 3,600.1 \\
\hline
\end{tabular}}
\caption{Solutions found and computational cost of solving the
  large-sized instances with learning rate $\alpha = 0.2$ using CPLEX
  and CP Optimizer with no warm start.}
\label{tab9}
\end{table}

\begin{table}[ht!]
\centering
\resizebox{\textwidth}{!}{
\rowcolors{5}{}{lightgray}
\begin{tabular}{cccccccccc}
\hline
\multirow{3}{*}{instance} 
&  & \multicolumn{4}{c}{IBM ILOG CPLEX} 
&  & \multicolumn{3}{c}{IBM ILOG CP Optimizer} \\
\cline{3-6} \cline{8-10} 
&  & \multirow{2}{*}{makespan} & \multicolumn{3}{c}{Effort measurement}
&  & \multirow{2}{*}{makespan} & \multicolumn{2}{c}{Effort measurement} \\
\cline{4-6} \cline{9-10} 
&  &  & \#iterations & \#B\&B Nodes & CPU &  &  & \#branches & CPU \\
\cline{1-1} \cline{3-6} \cline{8-10} 
DAFJS01 &  & {[}15,541, 24,191{]} 35.76\% & 3,223,970 & 12,860 & 3,600.0 &  & {[}10,711, 20,080{]} 46.66\% & 7,023,413 & 3,600.1 \\
DAFJS02 &  & {[}16,850, 27,398{]} 38.50\% & 3,318,330 & 12,085 & 3,600.0 &  & {[}10,500, 22,335{]} 52.99\% & 9,620,254 & 3,600.0 \\
DAFJS03 &  & {[}22,536, 165,343{]} 86.37\% & 125,755 & 7 & 3,600.0 &  & {[}21,152, 41,224{]} 48.69\% & 920,372 & 3,600.2 \\
DAFJS04 &  & {[}24,891, 120,636{]} 79.37\% & 197,181 & 48 & 3,600.0 &  & {[}23,278, 41,388{]} 43.76\% & 785,196 & 3,600.4 \\
DAFJS05 &  & {[}20,613, 101,212{]} 79.63\% & 482,040 & 586 & 3,600.0 &  & {[}12,401, 28,360{]} 56.27\% & 2,582,449 & 3,600.1 \\
\hline
\end{tabular}}
\caption{Solutions found and computational cost of solving the
  large-sized instances with learning rate $\alpha = 0.3$ using CPLEX
  and CP Optimizer with no warm start.}
\label{tab10}
\end{table}

\begin{table}[ht!]
\centering
\resizebox{0.8\textwidth}{!}{
\rowcolors{7}{}{lightgray}
\begin{tabular}{cccccccccc}
\hline
\multirow{3}{*}{instance}
&  & \multicolumn{4}{c}{IBM ILOG CPLEX}
&  & \multicolumn{3}{c}{IBM ILOG CP Optimizer} \\
\cline{3-6} \cline{8-10} 
&  & \multirow{2}{*}{makespan} & \multicolumn{3}{c}{Effort measurement}
&  & \multirow{2}{*}{makespan} & \multicolumn{2}{c}{Effort measurement} \\
\cline{4-6} \cline{9-10} 
&  &  & \#iterations & \#B\&B Nodes & CPU &  &  & \#branches & CPU \\
\cline{1-1} \cline{3-6} \cline{8-10} 
miniDAFJS01 &  & 22,875 & 27,130 & 272 & 2.2 &  & 22,875 & 143,152 & 3.2 \\
miniDAFJS02 &  & 22,708 & 9,960 & 85 & 1.1 &  & 22,708 & 15,317 & 0.3 \\
miniDAFJS03 &  & 18,363 & 0 & 0 & 0.0 &  & 18,363 & 5,403 & 0.1 \\
miniDAFJS04 &  & 20,498 & 155 & 0 & 0.0 &  & 20,498 & 27,577 & 0.2 \\
miniDAFJS05 &  & 20,593 & 27,452 & 195 & 1.7 &  & 20,593 & 66,727 & 1.1 \\
\hline
\end{tabular}}
\caption{Solutions found and computational cost of solving the
  small-sized instances with learning rate $\alpha = 0.1$ using CPLEX
  and CP Optimizer with warm start.}
\label{tab11}
\end{table}

\begin{table}[ht!]
\centering
\resizebox{0.85\textwidth}{!}{
\rowcolors{7}{}{lightgray}
\begin{tabular}{cccccccccc}
\hline
\multirow{3}{*}{instance}
&  & \multicolumn{4}{c}{IBM ILOG CPLEX}
&  & \multicolumn{3}{c}{IBM ILOG CP Optimizer} \\
\cline{3-6} \cline{8-10} 
&  & \multirow{2}{*}{makespan} & \multicolumn{3}{c}{Effort measurement}
&  & \multirow{2}{*}{makespan} & \multicolumn{2}{c}{Effort measurement} \\
\cline{4-6} \cline{9-10} 
&  &  & \#iterations & \#B\&B Nodes & CPU &  &  & \#branches & CPU \\
\cline{1-1} \cline{3-6} \cline{8-10} 
miniDAFJS01 &  & 21,327 & 539,475 & 4,331 & 31.9 &  & 21,327 & 759,456 & 13.6 \\
miniDAFJS02  &  & 20,635 & 30,346 & 207 & 2.0 &  & 20,635 & 30,019 & 0.5 \\
miniDAFJS03  &  & 17,972 & 561 & 0 & 0.1 &  & 17,972 & 10,221 & 0.3 \\
miniDAFJS04  &  & 19,602 & 655 & 0 & 0.1 &  & 19,602 & 14,122 & 0.2 \\
miniDAFJS05  &  & 18,803 & 58,235 & 514 & 3.5 &  & 18,803 & 206,620 & 4.5 \\
\hline
\end{tabular}}
\caption{Solutions found and computational cost of solving the
  small-sized instances with learning rate $\alpha = 0.2$ using CPLEX
  and CP Optimizer with warm start.}
\label{tab12}
\end{table}

\begin{table}[ht!]
\centering
\resizebox{0.85\textwidth}{!}{
\rowcolors{7}{}{lightgray}
\begin{tabular}{cccccccccc}
\hline
\multirow{3}{*}{instance}
&  & \multicolumn{4}{c}{IBM ILOG CPLEX}
&  & \multicolumn{3}{c}{IBM ILOG CP Optimizer} \\
\cline{3-6} \cline{8-10} 
&  & \multirow{2}{*}{makespan} & \multicolumn{3}{c}{Effort measurement}
&  & \multirow{2}{*}{makespan} & \multicolumn{2}{c}{Effort measurement} \\
\cline{4-6} \cline{9-10} 
&  &  & \#iterations & \#B\&B Nodes & CPU &  &  & \#branches & CPU \\
\cline{1-1} \cline{3-6} \cline{8-10} 
miniDAFJS01 &  & 19,443 & 698,096 & 4,945 & 38.9 &  & 19,443 & 1,528,271 & 32.0 \\
miniDAFJS02  &  & 18,916 & 43,060 & 324 & 2.0 &  & 18,916 & 81,435 & 1.5 \\
miniDAFJS03  &  & 17,419 & 6,368 & 91 & 0.6 &  & 17,419 & 27,973 & 0.4 \\
miniDAFJS04  &  & 18,800 & 4,986 & 135 & 0.3 &  & 18,800 & 23,757 & 0.3 \\
miniDAFJS05  &  & 17,596 & 174,306 & 2,110 & 7.9 &  & 17,596 & 505,368 & 9.0 \\
\hline
\end{tabular}}
\caption{Solutions found and computational cost of solving the
  small-sized instances with learning rate $\alpha = 0.3$ using CPLEX
  and CP Optimizer with warm start.}
\label{tab13}
\end{table}

\begin{table}[ht!]
\centering
\resizebox{\textwidth}{!}{
\rowcolors{5}{}{lightgray}
\begin{tabular}{cccccccccc}
\hline
\multirow{3}{*}{instance}
&  & \multicolumn{4}{c}{IBM ILOG CPLEX}
&  & \multicolumn{3}{c}{IBM ILOG CP Optimizer} \\
\cline{3-6} \cline{8-10} 
&  & \multirow{2}{*}{makespan} & \multicolumn{3}{c}{Effort measurement}
&  & \multirow{2}{*}{makespan} & \multicolumn{2}{c}{Effort measurement} \\
\cline{4-6} \cline{9-10} 
&  &  & \#iterations & \#B\&B Nodes & CPU &  &  & \#branches & CPU \\
\cline{1-1} \cline{3-6} \cline{8-10} 
DAFJS01 & & {[}21,618,   28,107{]} 23.09\% & 2,625,310 & 8,412 & 3,600.0 & & {[}18,588,   23,460{]} 20.77\% & 11,970,328 & 3,600.0 \\
DAFJS02 & & {[}23,418,   29,389{]} 20.32\% & 3,476,653 & 10,932 & 3,600.0 & & {[}18,533,   26,535{]} 30.16\% & 9,040,577 & 3,600.0 \\
DAFJS03 & & {[}42,594,   53,688$^\dagger${]} 20.66\% & 116,124 & 0 & 3,600.0 & & {[}41,271,   53,610{]} 23.02\% & 896,736 & 3,600.0 \\
DAFJS04 & & {[}45,501,   54,082$^\dagger${]} 15.87\% & 204,366 & 82 & 3,600.0 & & {[}44,062,   53,010{]} 16.88\% & 1,705,142 & 3,600.0 \\
DAFJS05 & & {[}30,528,   45,862$^\dagger${]} 33.44\% & 354,131 & 576 & 3,600.0 & & {[}23,378,   37,770{]} 38.1\% & 5,032,131 & 3,600.0 \\
\hline
\end{tabular}}
\caption{Solutions found and computational cost of solving the
  large-sized instances with learning rate $\alpha = 0.1$ using CPLEX
  and CP Optimizer with warm start.}
\label{tab14}
\end{table}

\begin{table}[ht!]
\centering
\resizebox{\textwidth}{!}{
\rowcolors{5}{}{lightgray}
\begin{tabular}{cccccccccc}
\hline
\multirow{3}{*}{instance}
&  & \multicolumn{4}{c}{IBM ILOG CPLEX}
&  & \multicolumn{3}{c}{IBM ILOG CP Optimizer} \\
\cline{3-6} \cline{8-10} 
&  & \multirow{2}{*}{makespan} & \multicolumn{3}{c}{Effort measurement}
&  & \multirow{2}{*}{makespan} & \multicolumn{2}{c}{Effort measurement} \\
\cline{4-6} \cline{9-10} 
&  &  & \#iterations & \#B\&B Nodes & CPU &  &  & \#branches & CPU \\
\cline{1-1} \cline{3-6} \cline{8-10} 
DAFJS01 & & {[}16,695,   26,067{]} 35.95\% & 1,867,058 & 5,324 & 3,600.0 & & {[}14,141,   21,861{]} 35.31\% & 12,852,459 & 3,600.0 \\
DAFJS02 & & {[}20,341,   28,154$^\dagger${]} 27.75\% & 577,583 & 1,370 & 3,600.0 & & {[}13,966,   24,746{]} 43.56\% & 18,869,445 & 3,600.0 \\
DAFJS03 & & {[}31,052,   48,555$^\dagger${]} 36.05\% & 84,156 & 0 & 3,600.0 & & {[}29,519,   46,422{]} 36.41\% & 1,344,012 & 3,600.0 \\
DAFJS04 & & {[}33,568,   48,461$^\dagger${]} 30.73\% & 183,960 & 44 & 3,600.0 & & {[}32,038,   45,677{]} 29.86\% & 1,406,727 & 3,600.0 \\
DAFJS05 & & {[}25,167,   44,790$^\dagger${]} 43.81\% & 324,450 & 373 & 3,600.0 & & {[}17,025,   33,314{]} 48.9\% & 3,620,699 & 3,600.0 \\
\hline
\end{tabular}}
\caption{Solutions found and computational cost of solving the
  large-sized instances with learning rate $\alpha = 0.2$ using CPLEX
  and CP Optimizer with warm start.}
\label{tab15}
\end{table}

\begin{table}[ht!]
\centering
\resizebox{\textwidth}{!}{
\rowcolors{5}{}{lightgray}
\begin{tabular}{cccccccccc}
\hline
\multirow{3}{*}{instance}
&  & \multicolumn{4}{c}{IBM ILOG CPLEX}
&  & \multicolumn{3}{c}{IBM ILOG CP Optimizer} \\
\cline{3-6} \cline{8-10} 
&  & \multirow{2}{*}{makespan} & \multicolumn{3}{c}{Effort measurement}
&  & \multirow{2}{*}{makespan} & \multicolumn{2}{c}{Effort measurement} \\
\cline{4-6} \cline{9-10} 
&  &  & \#iterations & \#B\&B Nodes & CPU &  &  & \#branches & CPU \\
\cline{1-1} \cline{3-6} \cline{8-10} 
DAFJS01 & & {[}13,085,   22,616$^\dagger${]} 42.14\% & 1,117,444 & 1,730 & 3,600.0 & & {[}10,730,   20,247{]} 47\% & 7,775,777 & 3,600.0 \\
DAFJS02 & & {[}16,801,   26,183$^\dagger${]} 35.83\% & 1,515,028 & 3,091 & 3,600.0 & & {[}10,500,   22,147{]} 52.59\% & 18,313,536 & 3,600.0 \\
DAFJS03 & & {[}22,719,   43,964$^\dagger${]} 48.32\% & 87,465 & 6 & 3,600.0 & & {[}21,152,   42,677{]} 50.44\% & 653,668 & 3,600.0 \\
DAFJS04 & & {[}24,927,   43,281$^\dagger${]} 42.41\% & 150,187 & 31 & 3,600.0 & & {[}23,278,   40,023{]} 41.84\% & 2,127,914 & 3,600.0 \\
DAFJS05 & & {[}20,474,   38,237$^\dagger${]} 46.46\% & 298,694 & 514 & 3,600.0 & & {[}12,401,   29,897{]} 58.52\% & 3,366,879 & 3,600.0 \\
\hline
\end{tabular}}
\caption{Solutions found and computational cost of solving the
  large-sized instances with learning rate $\alpha = 0.3$ using CPLEX
  and CP Optimizer with warm start.}
\label{tab16}
\end{table}

\section{Conclusions}

In this work, we addressed the FJS scheduling problem with sequencing
flexibility and position-based learning effect. To the authors'
knowledge, this combination, with potential application in a wide
range of real-world industrial environments, has never been addressed
in the literature. As a first step towards its efficient and effective
solution, we introduced a set of 110 instances that transform into 330
instances by varying the learning rate $\alpha \in \{0.1, 0.2,
0.3\}$. By introducing MILP and CP models, an exact solver aided by
constructive heuristics was able to provide 183 proven optimal
solutions. Instances, their solutions, models and constructive
heuristics are all available for download at
\url{https://github.com/kennedy94/FJS}. We expect this benchmark
test-set to guide the introduction and evaluation of new effective
heuristic and metaheuristic approaches for the considered problem. In
fact, this is the current line of research of the authors.\\

\noindent
\textbf{Conflict of interest statement:} On behalf of all authors, the
corresponding author states that there is no conflict of interest.\\

\noindent
\textbf{Data availability:} The datasets generated during and/or
analysed during the current study are available in the GitHub
repository, \url{https://github.com/kennedy94/FJS}.

\bibliographystyle{plain}

\bibliography{arbro2023}

\end{document}